\documentclass[11pt, a4paper]{article}
\usepackage{latexsym,amssymb}
\usepackage{amsmath}
\usepackage[mathscr]{eucal}
\usepackage{color, xcolor}
\usepackage[abbrev]{amsrefs}

\newtheorem{Theorem}{\bf Theorem}[section]
\newtheorem{Lemma}{\bf Lemma}[section]
\newtheorem{Proposition}{\bf Proposition}[section]
\newtheorem{Corollary}{\bf Corollary}[section]
\newtheorem{Remark}{\bf Remark}[section]
\newtheorem{Example}{\bf Example}[section]
\newtheorem{Definition}{\bf Definition}[section]

\newenvironment{theorem}{\begin{Theorem}$\!\!\!$}{\end{Theorem}}
\newenvironment{lemma}{\begin{Lemma}$\!\!\!$}{\end{Lemma}}
\newenvironment{proposition}{\begin{Proposition}$\!\!\!$}{\end{Proposition}}

\newenvironment{remark}{\begin{Remark}$\!\!\!$}{\end{Remark}}

\newenvironment{definition}{\begin{Definition}$\!\!\!$}{\end{Definition}}

\numberwithin{equation}{section}

\newcommand{\dee}{{\rm{d}}}
\newcommand{\loc}{{\rm{loc}}}
\newcommand{\R}{\mathbb{R}}
\newcommand{\N}{\mathbb{N}}
\newcommand{\G}{\mathbb{G}_{\Omega}}

\newcommand{\EE}{\mathcal{E}}

\newcommand{\dt}{\dfrac{\dee}{\dee{\mathnormal t}}}

\newcommand{\dtau}{\dfrac{\dee}{\dee\tau}}

\newcommand*{\RL}{\partial_{\mathnormal t}^{\alpha}}

\DeclareMathOperator*{\esssup}{ess\sup}

\begin{document}
\title{Large-Time Behavior and Asymptotic Profiles for \\
Time-Fractional Nonlinear Diffusion Equations \\
via Entropy Methods}
\author{
Mikiya Kametaka \\
\\
Graduate School of Applied Mathematics and Informatics Course,\\ 
Faculty of Advanced Science and Technology, Ryukoku University,\\
1-5 Yokotani, Seta Oe-cho, Otsu, Shiga 520-2194, Japan.\\
\\
E-mail: {\tt Y25D001@mail.ryukoku.ac.jp}
}
\date{}
\maketitle
\begin{abstract}
We investigate the large-time behavior of nonnegative solutions to the time-fractional nonlinear diffusion equation
$\RL(u-u_0) = \Delta u^m$, $0<\alpha<1$, $m>0$,
in a bounded domain with the homogeneous Dirichlet boundary condition.
We prove that the rescaled solution $t^{\frac{\alpha}{m}}u$ converges as $t\to\infty$
to an asymptotic profile $S$ characterized by $-\Delta S^m = u_0/\Gamma(1-\alpha)$.
Moreover, we show that $t^\alpha u^m$ converges to $S^m$, which admits an explicit representation in terms of the Green function and the initial datum.
In contrast to the classical case $\alpha=1$, the asymptotic profile depends explicitly on the initial datum.
\end{abstract}
\vspace{10pt}

\vspace{40pt}
\noindent 



\newpage

\section{Introduction}
\label{Section:1}
This paper is concerned with the large-time behavior of nonnegative solutions to the initial-boundary value problem
\begin{equation}
\label{P}\tag{P}
	\left\{
	\begin{array}{lll}
	\RL(u-u_0) = \Delta u^m,
	&
	x \in \Omega, \quad t>0,
	\vspace{5pt}\\
	u = 0,
	&
	x \in \partial\Omega, \quad t>0,
	\vspace{5pt}\\
	u(x,0) = u_0(x)\ge 0,
	&
	x \in \Omega,
	\end{array}
	\right.
\end{equation}
where $m>0$, $0<\alpha<1$, $N \ge 3$, and $\Omega\subset\R^N$ is a bounded domain with smooth boundary.
Moreover, $\RL$ stands for the $\alpha$-th order Riemann-Liouville derivative, which is defined by
$$
	\RL v(t):=\dt\bigl(g_{1-\alpha}*v\bigr)(t),
	\qquad
	g_{1-\alpha}(t):=\frac{t^{-\alpha}}{\Gamma(1-\alpha)}.
$$
Here $*$ denotes the convolution on $(0,\infty)$ with respect to the time variable, that is,
$$
	(f*g)(t):=\int_0^t f(t-\tau)g(\tau)\,\dee\tau,
	\qquad t\ge0
$$
for $f,g\in L^1_{\rm loc}([0,\infty))$.

Throughout this paper, we use the following notion of an $H^{-1}$-solution introduced in \cite{A,BGI}.
\begin{definition}
\label{Definition:1.1}
	We say that $u \in L_{\loc}^2(0,\infty; H^{-1}(\Omega))$ is an $H^{-1}$-solution to problem \eqref{P} if
	\begin{itemize}
	\item[$({\rm i})$]
	$g_{1-\alpha}*[u-u_0] \in H_{\loc}^1(0,\infty; H^{-1}(\Omega))$ with $\bigl(g_{1-\alpha}*[u-u_0]\bigr)(0)=0$,
	\item[$({\rm ii})$]
	$u^m \in L_{\loc}^2(0,\infty; H_0^1(\Omega))$,
	\item[$({\rm iii})$]
	$u$ satisfies the identity
	\begin{equation}
	\label{eq:1.1}
		\int_{\Omega}\RL(u(x,\cdot)-u_0(x))(t)\varphi(x)\,\dee x
		+ \int_{\Omega}\nabla u^m(x,t)\cdot\nabla\varphi(x)\,\dee x = 0,
		\quad \mbox{a.e. $t>0$}
	\end{equation}
	for any test function $\varphi\in H_0^1(\Omega)$.
	\end{itemize}
\end{definition}
\bigskip

The existence of solutions to time-fractional nonlinear diffusion equations has been studied extensively (see, e.g., \cite{A,AN,BGI,KK,WWZ}).
In particular, Wittbold, Wolejko, and Zacher \cite{WWZ} established the existence of bounded weak solutions for bounded initial data in the degenerate case $m>1$.
The author and Kawakami \cite{KK} extended this result to the singular case $0<m<1$ and gave improved regularity properties of weak solutions.
Furthermore, Akagi \cite{A} developed the abstract theory for time-fractional gradient flows of subdifferential type in a real Hilbert space $H$
$$
	\RL(u-u_0)(t)+\partial\varphi(u(t))\ni0
	\quad \mbox{in $H$},\quad t>0,
$$
and proved the existence and uniqueness of strong solutions (see also \cite{AN}).
Here $\partial\varphi$ denotes the subdifferential operator of $\varphi$.
We can apply this framework to problem~\eqref{P}.
By the results in \cite{A, BGI, LL}, for any nonnegative initial datum $u_0 \in L^{m+1}(\Omega) \cap H^{-1}(\Omega)$, there exists a unique nonnegative $H^{-1}$-solution
$u \in L^\infty(0,\infty; L^{m+1}(\Omega) \cap H^{-1}(\Omega)) \cap C([0,\infty); H^{-1}(\Omega))$ to problem~\eqref{P}.
Moreover, it holds that
\begin{equation}
\label{eq:1.2}
	\|u(t)\|_{m+1} \le \|u_0\|_{m+1}, \quad \mbox{a.e. $t>0$},
\end{equation}
where $\|\cdot\|_p = \|\cdot\|_{L^p(\Omega)}$ with $1\le p \le \infty$.

We now turn to the large-time behavior.
Time decay estimates for time-fractional diffusion equations have been established by many authors (see, e.g., \cite{BT,CQW,KSVZ,KSZ,VZ2,ST}).
In particular, Vergara and Zacher \cite{VZ2} showed that the $L^2(\Omega)$-norm of weak solutions to linear equations
with measurable coefficients
$$
	\RL(u-u_0) - {\rm div}(A(x,t)\nabla u) = 0
$$
decays like $t^{-\alpha}$.
Furthermore, Bonforte, Gualdani, and Ibarrondo \cite{BGI} proved the nonextinction property of nonnegative solutions to problem~\eqref{P} and
established the asymptotic behavior
$$
	\|u(t)\|_p \asymp t^{-\frac{\alpha}{m}}
	\quad \mbox{as $t\to\infty$}
$$
for every $p \in [1,\infty]$ (see Proposition~\ref{Proposition:4.1} below for details).
For the linear equation $m=1$ in the whole space, namely,
$$
	\RL(u-u_0) = \Delta u, \quad x\in\R^N,\,\, t>0, \qquad u(x,0) = u_0(x), \quad x\in\R^N,
$$
Cort\'azar, Quir\'os, and Wolanski~\cite{CQW} showed that the rescaled solution $t^\alpha u(x,t)$ converges to
the Newtonian potential of the initial datum on compact sets in $\R^N$.

In the classical case, namely, $\alpha=1$,
for the singular case $0<m<1$, solutions extinguish in finite time,
while for the degenerate case $m>1$ they decay polynomially and their sharp rates are well known (see, e.g., survey books \cite{V1,V2}).
After a suitable rescaling, problem~\eqref{P} is transformed into
$$
	\left\{
	\begin{array}{lll}
	\partial_\tau v = \Delta v^m + \lambda v,
	&
	x\in\Omega, \quad \tau>0,
	\vspace{5pt}\\
	v = 0,
	&
	x\in\partial\Omega, \quad \tau>0,
	\vspace{5pt}\\
	v(x,0) = v_0(x),
	&
	x\in\Omega.
	\end{array}
	\right.
$$
Based on the theory developed by Feireisl and Simondon \cite{FS},
for the fast diffusion range $m_s<m<1$,
Bonforte, Grillo, and V\'azquez \cite{BGV} proved the convergence of the rescaled solution $v$ to a positive stationary profile $S$ satisfying
$$
	-\Delta S^m=\lambda S \quad \mbox{in $\Omega$},
	\qquad
	S=0 \quad \mbox{on $\partial\Omega$},
$$
where $m_s = [N-2]_+/(N+2)$ is the inverse Sobolev exponent.
Furthermore, V\'azquez \cite{V3} gave the corresponding result for the porous medium range $m>1$ (see also \cite{AP}).
Moreover, related asymptotic results for the fast diffusion equation of the form
$$
	\partial_t(|u|^{q-2}u) = \Delta u, \quad q>2
$$
have also been obtained by energy methods (see, e.g., \cite{A1,AM,AK} and the references therein).
Thus, in the classical case, suitably rescaled solutions converge to the spatial profile of a separable solution.

On the other hand, the time-fractional problem exhibits significantly different properties.
Due to the time nonlocality of fractional derivatives,
the standard chain and product rules, as well as the usual rescaling arguments, are not directly applicable.
For instance, for sufficiently smooth functions $u$ and $\varphi$, we have
$$
	\varphi(t)\RL u(t) = \RL(\varphi u)(t)
	- \frac{\alpha}{\Gamma(1-\alpha)}\int_0^t\frac{\varphi(t)-\varphi(s)}{(t-s)^{\alpha+1}}u(s)\,\dee s.
$$
This identity follows from the corresponding product rule for nonlocal derivatives (see Lemma~\ref{Lemma:2.2} below for details).
The additional nonlocal term illustrates the difficulty of applying time-dependent rescaling arguments to time-fractional diffusion equations.
In particular, the identification of the asymptotic profile for the time-fractional nonlinear problem~\eqref{P} remains an open problem.
The aim of the present paper is to identify a profile $S$ for problem~\eqref{P} and to show that,
in contrast to the classical case, $S$ depends explicitly on the initial datum.
We now state our main result.
\begin{theorem}
\label{Theorem:1.1}
	Assume that $0<\alpha<1$, $m>0$, $N\ge3$, and $p_0>1$ satisfies
	$$
		\left\{
		\begin{array}{ll}
		p_0 \ge m+1
		&
		\mbox{if} \quad m+1>\dfrac{N}{2},
		\vspace{5pt}\\
		p_0 > \dfrac{N}{2}
		&
		\mbox{if} \quad \dfrac{N}{2}\ge m+1.
		\end{array}
		\right.
	$$
	Let $u \in L^\infty(0,\infty; L^{m+1}(\Omega) \cap H^{-1}(\Omega)) \cap C([0,\infty); H^{-1}(\Omega))$ be a nonnegative $H^{-1}$-solution to problem~\eqref{P}
	with nonnegative initial datum $u_0 \in L^{p_0}(\Omega) \cap H^{-1}(\Omega)$.
	Moreover, let $S$ be the unique nonnegative weak solution to problem
	\begin{equation}
	\label{eq:1.3}
		-\Delta S^m = c_{\alpha}u_0 \quad
		\mbox{in $\Omega$},
		\qquad
		S = 0 \quad
		\mbox{on $\partial\Omega$},
	\end{equation}
	where $c_{\alpha} := 1/\Gamma(1-\alpha)$.
	Then, for every $1 \le p < \infty$,
	\begin{equation}
	\label{eq:1.4}
		\lim_{t\to\infty}\|t^{\frac{\alpha}{m}}u(t) - S\|_p = 0.
	\end{equation}
	Furthermore, it holds that
	\begin{equation}
	\label{eq:1.5}
		\lim_{t\to\infty}\left\|t^{\alpha}u^m(t) - c_{\alpha}\int_{\Omega}\G(\cdot,y)u_0(y)\,\dee y\right\|_p = 0
	\end{equation}
	for every $1 \le p < \infty$, where $\G$ is the Green function of $-\Delta_{\Omega}$.
\end{theorem}

\begin{remark}
\label{Remark:1.1}
	\begin{itemize}
	\item[$({\rm i})$]
	Since $u_0 \in H^{-1}(\Omega)$, the Lax-Milgram theorem guarantees the existence and uniqueness of weak solutions to problem
	$$
		-\Delta V = c_{\alpha}u_0 \quad
		\mbox{in $\Omega$},
		\qquad
		V = 0 \quad
		\mbox{on $\partial\Omega$}.
	$$
	Moreover, since $u_0 \ge 0$, by the maximum principle, it follows that $V \ge 0$,
	so we may define $S^m := V$ as a weak solution to \eqref{eq:1.3} satisfying
	\begin{equation}
	\label{eq:1.6}
		\int_{\Omega}\nabla S^m(x)\cdot\nabla\varphi(x)\,\dee x
		= c_{\alpha}\int_{\Omega}u_0(x)\varphi(x)\,\dee x
	\end{equation}
	for any test function $\varphi \in H_0^1(\Omega)$ $($see, e.g., \cite{GT}$)$.
	
	\item[$({\rm ii})$]
	Since $u_0$ also belongs to $L^{p_0}(\Omega)$ with $p_0 > N/2$,
	$S^m$ admits the representation
	$$
		S^m(x) = c_{\alpha}\int_{\Omega}\G(x,y)u_0(y)\,\dee y.
	$$
	Then it holds that $\|S\|_\infty \le C\|u_0\|_{p_0}^{1/m}$,
	where $C>0$ depends only on $\alpha$, $m$, $N$, and $\Omega$ $($see, e.g., \cite{BGI}*{Lemma~4.5}$)$.
	
	\item[$({\rm iii})$]
	When $m=1$, the convergence \eqref{eq:1.4}, namely, \eqref{eq:1.5} can also be proved directly from the representation formula of solutions
	$($see Section~$\ref{Section:5}$ below for details$)$.
	\end{itemize}
\end{remark}
\bigskip

We briefly explain the idea of the proof of Theorem~\ref{Theorem:1.1}.
The proof is based on a nonlinear entropy method introduced in \cite{BF}.
More precisely, we set $v(x,t):=t^{\frac{\alpha}{m}}u(x,t)$ and introduce the entropy functional
$$
	\EE[v|S](t)
	:= \int_\Omega\left\{\frac{v^{m+1}(x,t)}{m+1}-\frac{S^{m+1}(x)}{m+1}-S^m(x)(v(x,t)-S(x))\right\}\,\dee x.
$$
Using the Yosida approximation of the nonlocal time-derivative operator and a nonlocal chain rule inequality,
we derive a fractional differential inequality for $\EE[v|S]$.
Next, we compare the entropy functional with solutions to fractional ODEs and obtain decay estimates for the entropy functional.
Finally, by the coercivity of the entropy functional, we can derive the convergence results \eqref{eq:1.4} and \eqref{eq:1.5}.
In view of the Green function representation of $S^m$, the convergence \eqref{eq:1.5} implies that the asymptotic behavior can be expressed in terms of the initial datum.
This is consistent with the asymptotic behavior obtained in \cite{CQW} for the linear problem in the whole space.
\medskip

The rest of the present paper is organized as follows.
In Section~\ref{Section:2}, we collect the tools needed in the proof of the main result.
In particular, we recall properties of nonlocal time derivatives.
In Section~\ref{Section:3}, we establish decay estimates for fractional ODEs.
In Section~\ref{Section:4}, we prove Theorem~\ref{Theorem:1.1}.
We first derive a basic entropy inequality for the rescaled solution.
Due to the time nonlocality, this inequality exhibits different behavior depending on whether $m>\alpha$, $m=\alpha$, or $0<m<\alpha$.
Therefore, we divide the proof of Theorem~\ref{Theorem:1.1} into four cases:
$m\ge1$, $\alpha<m<1$, $m=\alpha$, and $0<m<\alpha$.
Finally, in Section~\ref{Section:5}, we give an alternative proof in the linear case $m=1$ based on the representation formula of solutions.

\section{Preliminaries}
\label{Section:2}
In this section, we first recall some properties of the Mittag-Leffler function, which is defined by
$$
	E_{a,b}(z) = \sum_{k=0}^\infty\frac{z^k}{\Gamma(a k+b)}, \quad z\in\R, \quad a,b>0.
$$
This function satisfies
\begin{equation}
\label{eq:2.1}
	\frac{1}{1+\Gamma(1-\alpha)t} \le E_{\alpha,1}(-t) \le \frac{1}{1+\frac{1}{\Gamma(1+\alpha)}t},
	\quad t>0
\end{equation}
for $0<\alpha<1$.
Moreover, it holds that
\begin{equation}
\label{eq:2.2}
	\dt E_{\alpha,1}(-\lambda t^{\alpha}) = -\lambda t^{\alpha-1}E_{\alpha,\alpha}(-\lambda t^{\alpha}),
	\quad t>0
\end{equation}
as well as that
\begin{equation}
\label{eq:2.3}
	\frac{1}{(1+c_{\alpha,1}t)^2} \le \Gamma(\alpha)E_{\alpha, \alpha}(-t) \le \frac{1}{(1+c_{\alpha,2}t)^2},
	\quad t\ge0,
\end{equation}
where $c_{\alpha,1}, c_{\alpha,2}$ are positive constants depending only on $\alpha$ (see, e.g., \cite{BS}).

We next recall some basic facts of the time-fractional derivative operator $\RL$.
Let $T>0$, $1\le p<\infty$, and let $X$ be a real Banach space.
We define the operator $B$ by
$$
	\begin{aligned}
	&
	Bu = \dt(g_{1-\alpha}*u), \\
	&
	D(B) = \left\{u \in L^p(0,T; X);\, g_{1-\alpha}*u \in W^{1,p}(0,T; X)\ \mbox{with}\ (g_{1-\alpha}*u)(0)=0\right\}.
	\end{aligned}
$$
It is well known that $B$ is an $m$-accretive operator in $L^p(0,T; X)$ (see \cites{B,GLS,G}).
Let $n \in \N$, and let $B_n$ denote the Yosida approximation of $B$, which is defined by $B_n = nB(n+B)^{-1}$.
Then, for any $u \in D(B)$, it holds that
\begin{equation}
\label{eq:2.4}
	B_nu \to Bu
	\quad\mbox{strongly in $L^p(0,T; X)$ as $n \to \infty$}.
\end{equation}
Moreover, $B_n$ admits the representation
$$
	B_nu = \dt(g_{1-\alpha,n}*u), \quad u \in L^p(0,T;X),
$$
where $g_{1-\alpha,n} = ns_{\alpha,n}$, and $s_{\alpha,n}$ is the unique solution to the Cauchy problem
$$
	\dt\bigl(g_{1-\alpha}*[s_{\alpha,n}-s_{\alpha,n}(0)]\bigr) + ns_{\alpha,n} = 0, \quad t>0,
	\qquad s_{\alpha,n}(0)=1.
$$
See, e.g., \cite{VZ1}.
This problem can also be written as the Volterra integral equation
\begin{equation}
\label{eq:2.5}
	s_{\alpha,n}(t) + n(s_{\alpha,n}*g_{\alpha})(t) =1,\quad t>0,
\end{equation}
and, using the Laplace transform, the solution $s_{\alpha,n}$ is represented by
$$
	s_{\alpha,n}(t) = E_{\alpha,1}(-nt^{\alpha}), \quad t \ge 0.
$$
Therefore, by \eqref{eq:2.1}-\eqref{eq:2.3},
$s_{\alpha,n} \in W_{\loc}^{1,1}([0,\infty))$ is a nonnegative and nonincreasing function.
Let $h_{\alpha,n} \in L_{\rm loc}^1([0,\infty))$ be the nonnegative resolvent kernel associated with $ng_{\alpha}$, that is,
\begin{equation}
\label{eq:2.6}
	h_{\alpha,n}(t) + n(h_{\alpha,n}*g_{\alpha})(t) = ng_{\alpha}(t),\quad t>0
\end{equation}
(see, e.g., \cite{GLS}).
Then it follows from $g_{1-\alpha}*g_{\alpha}=1$ that
$$
	(g_{1-\alpha}*h_{\alpha,n})(t) + n(g_{1-\alpha}*h_{\alpha,n}*g_{\alpha})(t) = n,\quad t>0,
$$
which together with \eqref{eq:2.5} implies that
$$
	g_{1-\alpha,n} = ns_{\alpha,n} = g_{1-\alpha}*h_{\alpha,n}.
$$
Note that, for any function $f \in L^p(0,T; X)$, it holds that
\begin{equation}
\label{eq:2.7}
	h_{\alpha,n}*f \to f\quad\mbox{strongly in $L^p(0,T; X)$ as $n \to \infty$}.
\end{equation}
In fact, set $v:=g_{\alpha}*f$. 
Then, by \eqref{eq:2.4}, we see that $v \in D(B)$ and
$$
	\dt(g_{1-\alpha,n}*v)=B_nv 
	\xrightarrow[n\to\infty]{}
	Bv=\dt(g_{1-\alpha}*v) = \dt(g_{1-\alpha}*g_{\alpha}*f) = f\quad\mbox{in $L^p(0,T;X)$},
$$
while also
$$
	B_nv = \dt(g_{1-\alpha}*h_{\alpha,n}*g_{\alpha}*f) = \dt(g_{1-\alpha}*g_{\alpha}*h_{\alpha,n}*f) = h_{\alpha,n}*f.
$$
This implies that \eqref{eq:2.7} holds.
In particular, putting $f=g_{1-\alpha}$ and $X=\R$, we have
\begin{equation}
\label{eq:2.8}
	g_{1-\alpha,n} \to g_{1-\alpha}
	\quad\mbox{strongly in $L^1(0,T)$ as $n \to \infty$}.
\end{equation}

Next, we recall a chain rule inequality and a product rule for nonlocal time-derivative operators.
The following lemmas play a crucial role to give a basic inequality for the entropy functional.
\begin{lemma}
\label{Lemma:2.1}\cite{KSVZ}*{Lemma~6.1}
	Let $T >0$ and $I\subset\R$ be open.
	Let further $h \in W^{1,1}(0,T)$, 
	$H \in C^1(I)$, $u_0 \in I$, and $u \in L^1(0,T)$ with $u(t) \in I$ for a.e. $t \in (0,T)$.
	Suppose that the functions $H(u), H'(u)u$, and $H'(u)(h'*u)$ belong to $L^1(0,T)$.
	Then it holds that
	$$
	\begin{aligned}
		&
		H'(u(t))\dt\bigl(h*u(\cdot)\bigr)(t)
		= \dt\bigl(h*H(u(\cdot))\bigr)(t) + \biggl(-H(u(t))+H'(u(t))u(t)\biggr)h(t) \\
		&
		+ \int_0^t\biggl(H(u(t-s))-H(u(t))-H'(u(t))\bigl(u(t-s)-u(t)\bigr)\biggr)\bigl[-h'(s)\bigr]\,\dee s,
		\quad \mbox{a.e. $t\in(0,T)$}.
	\end{aligned}
	$$
	In particular, if $H \in C^1(\R)$ with $H(0)=0$ is convex and $h$ is nonnegative and nonincreasing,
	it holds that
	\begin{equation}
	\label{eq:2.9}
		H'(u(t))\dt\bigl(h*u(\cdot)\bigr)(t) \ge \dt\bigl(h*H(u(\cdot))\bigr)(t), \quad \mbox{a.e. $t\in(0,T)$}.
	\end{equation}
\end{lemma}
\medskip

\begin{lemma}
\label{Lemma:2.2}
\cite{Z1}*{Lemma~2.2}
	Let $T>0$, $h \in W^{1,1}(0,T)$, $u \in L^1(0,T)$, and $\varphi \in C^1([0,T])$.
	Then it holds that
	$$
		\varphi(t)\dt(h*u)(t) = \dt(h*[\varphi u])(t) + \int_0^t h'(t-s)\bigl(\varphi(t)-\varphi(s)\bigr)u(s)\,\dee s,
		\quad \mbox{a.e. $t\in(0,T)$}.
	$$
\end{lemma}
\begin{remark}
\label{Remark:2.1}
In Lemma~$\ref{Lemma:2.2}$, the assumption
$\varphi\in C^1([0,T])$ can be relaxed.
More precisely, if
$$
    h\in W^{1,1}(0,T),
    \quad
    \varphi\in L^\infty(0,T),
    \quad
    u \in L^1(0,T),
$$
then it follows that
\begin{equation}
\label{eq:2.10}
	\varphi(t)\dt(h*u)(t)
	=
	\dt\bigl(h*[\varphi u]\bigr)(t)
	+
	\int_0^th'(t-s)\bigl(\varphi(t)-\varphi(s)\bigr)u(s)\,\dee s
\end{equation}
for a.e. $t\in(0,T)$.
Indeed, since
$$
	h*f \in W^{1,1}(0,T), \qquad \dt\bigl(h*f\bigr)(t) = \bigl(h'*f\bigr)(t) + h(0)f(t), \quad \mbox{a.e. $t\in(0,T)$}
$$
for any $f \in L^1(0,T)$ and $h \in W^{1,1}(0,T)$,
we deduce that
$$
	\begin{aligned}
	&
	\varphi(t)\dt(h*u)(t) \\
	&
	= h(0)\varphi(t)u(t) + \varphi(t)(h'*u)(t) \\
	&
	= h(0)\varphi(t)u(t) + (h'*[\varphi u])(t) - \int_0^th'(t-s)\varphi(s)u(s)\,\dee s + \int_0^th'(t-s)\varphi(t)u(s)\,\dee s \\
	&
	= \dt(h*[\varphi u])(t) + \int_0^th'(t-s)\bigl(\varphi(t)-\varphi(s)\bigr)u(s)\,\dee s,
	\quad \mbox{a.e. $t\in(0,T)$}.
	\end{aligned}
$$
In particular, this extension is applicable to the regularized kernel $g_{1-\alpha,n}$.
\end{remark}
\bigskip

Finally, we recall some elementary inequalities for the following nonnegative convex function.
\begin{lemma}
\label{Lemma:2.3}
	Let $a,b \ge 0$ and $\theta > 0$.
	Define the nonnegative function $F_{\theta}$ by
	$$
		F_{\theta}(a,b) = \frac{1}{\theta+1}a^{\theta+1} - \frac{1}{\theta+1}b^{\theta+1} - b^\theta(a-b).
	$$
	Then the following assertions hold.
	\begin{itemize}
	\item[$({\rm i})$]
	If $\theta\ge1$, it follows that
	$$
		F_{\theta}(a,b) \ge \frac{1}{\theta+1}|a-b|^{\theta+1}.
	$$
	\item[$({\rm ii})$]
	If $0<\theta<1$, it follows that
	$$
		F_{\theta}(a,b) \le \frac{M^{1-\theta}}{2\theta}|a^\theta-b^\theta|^2, \quad M = \max\{a,b\}.
	$$
	\end{itemize}
\end{lemma}
{\bf Proof.}
We first prove assertion (i). Assume that $a\ge b$.
Since
$$
	x^\theta-y^\theta\ge (x-y)^\theta,
	\qquad x\ge y\ge0,\ \theta\ge1,
$$
we have
$$
	\begin{aligned}
	F_\theta(a,b)
	&
	=\int_b^a(r^\theta-b^\theta)\,\dee r
	\ge \int_b^a(r-b)^\theta\,\dee r
	=\frac{1}{\theta+1}(a-b)^{\theta+1}.
	\end{aligned}
$$
The case $a<b$ is analogous, and hence assertion (i) follows.

Next, let $0<\theta<1$, and again assume that $a\ge b$.
By the change of variables \(s=r^\theta\),
we obtain
$$
	F_\theta(a,b)
	=\frac1\theta\int_{b^\theta}^{a^\theta}
	(s-b^\theta)s^{\frac1\theta-1}\,\dee s
	\le \frac{a^{1-\theta}}{\theta}\int_{b^\theta}^{a^\theta}(s-b^\theta)\,\dee s
	= \frac{a^{1-\theta}}{2\theta}(a^\theta-b^\theta)^2.
$$
The case $a<b$ is treated in the same way with $a$ replaced by $b$.
Therefore, assertion (ii) holds.
This completes the proof.
$\Box$

\section{Decay estimates for fractional ODEs}
\label{Section:3}
In this section, we establish decay estimates for fractional ODEs.
We first recall an asymptotic decay estimate for a nonlinear fractional ODE with time-dependent coefficients.
This result will be applied to the case $m>\alpha$.
\begin{proposition}
\label{Proposition:3.1}
\cite{ST}*{Lemma~2.1}
	Let $0<\alpha<1$, $\beta > -\alpha$, $\gamma>0$, and $\lambda>0$.
	Assume that $w \in C([0,\infty))$ is the solution to the Cauchy problem
	$$
		\RL(w-w_0) + \lambda t^{\beta}w^\gamma = 0, \quad t>0,
		\qquad w(0) = w_0>0.
	$$
	Then it holds that
	$$
		\frac{c_0}{1+t^{\frac{\alpha+\beta}{\gamma}}} \le w(t) \le \frac{c_1}{1+t^{\frac{\alpha+\beta}{\gamma}}},
		\quad t>0,
	$$
	where $c_0$, $c_1$ are positive constants depending only on $\alpha$, $\beta$, $\gamma$, $\lambda$, and $w_0$.
\end{proposition}
\bigskip

In order to give decay estimates for fractional ODEs,
we first develop the comparison principle for inhomogeneous fractional ODEs with time-dependent coefficients.
By applying the same arguments as in \cite{VZ2}*{Lemma~2.6}, we obtain the following lemma.
\begin{lemma}
\label{Lemma:3.1}
	Let $f \in C(\R)$ be a nondecreasing function, and $r, g \in L_{\loc}^1([0,\infty)) \cap C(0,\infty)$,
	and let $r$ be a nonnegative function.
	Assume that $u, v \in C([0,\infty))$ satisfy
	$$
		\begin{array}{ll}
		\RL(u-u_0) + r(t)f(u) \le g, \quad t>0,
		&
		\quad u(0) = u_0,
		\vspace{5pt}\\
		\RL(v-v_0) + r(t)f(v) \ge g, \quad t>0,
		&
		\quad v(0)=v_0.
		\end{array}
	$$
	If $u_0 \le v_0$, then $u(t) \le v(t)$ for all $t>0$.
\end{lemma}
The proof of this lemma is very similar to \cite{VZ2}*{Lemma~2.6}.
However, for the reader's convenience, we write the details.

\noindent
{\bf Proof of Lemma~\ref{Lemma:3.1}.}
Let $T>0$ be arbitrary.
Convolving the inequality for $u$ with the nonnegative kernel $h_{\alpha,n}$, we have
\begin{equation}
\label{eq:3.1}
	\dt\bigl(g_{1-\alpha,n}*[u-u_0]\bigr)(t) + \biggl(h_{\alpha,n}*\biggl[r(\cdot)f(u(\cdot))\biggr]\biggr)(t) \le (h_{\alpha,n}*g)(t),
	\quad 0<t<T.
\end{equation}
Indeed, since $\bigl(g_{1-\alpha}*[u-u_0]\bigr)(0) = 0$, it holds that
$$
	\begin{aligned}
	\biggl(h_{\alpha,n}*\biggl[\dt\bigl(g_{1-\alpha}*[u-u_0]\bigr)\biggr]\biggr)(t)
	&
	= \dt\bigl(h_{\alpha,n}*g_{1-\alpha}*[u-u_0]\bigr)(t) - h_{\alpha,n}(t)\bigl(g_{1-\alpha}*[u-u_0]\bigr)(0) \\
	&
	= \dt\bigl(g_{1-\alpha,n}*[u-u_0]\bigr)(t),
	\quad 0<t<T.
	\end{aligned}
$$
In the same way, we convolve the inequality for $v$ with $h_{\alpha,n}$ and take the difference with \eqref{eq:3.1} to get
$$
	\dt\bigl(g_{1-\alpha,n}*[u-v]\bigr)(t) + \biggl(h_{\alpha,n}*\biggl[r(\cdot)\biggl(f(u(\cdot))-f(v(\cdot))\biggr)\biggr]\biggr)(t) \le g_{1-\alpha,n}(t)(u_0-v_0) \le 0
$$
for $0<t<T$.
Let $[y]_+$ be a positive part of $y \in \R$,
multiplying this inequality by $[u-v]_+$ and using \eqref{eq:2.9} in Lemma~\ref{Lemma:2.1} for $H(y) = \frac{1}{2}[y]_+^2$,
we deduce that
$$
	\frac{1}{2}\dt\bigl(g_{1-\alpha,n}*[u-v]_+^2\bigr)(t)
	+ \biggl(h_{\alpha,n}*\biggl[r(\cdot)\biggl(f(u(\cdot))-f(v(\cdot))\biggr)\biggr]\biggr)(t)[u(t)-v(t)]_+
	\le 0,
	\quad 0<t<T.
$$
Since $\bigl(g_{1-\alpha,n}*[u-v]_+^2\bigr)(0)=0$, integrating from $0$ to $T$, we obtain
\begin{equation}
\label{eq:3.2}
	\frac{1}{2}\bigl(g_{1-\alpha,n}*[u-v]_+^2\bigr)(T)
	+ \int_0^T\biggl(h_{\alpha,n}*\biggl[r(\cdot)\biggl(f(u(\cdot))-f(v(\cdot))\biggr)\biggr]\biggr)(t)[u(t)-v(t)]_+\,\dee t
	\le 0.
\end{equation}
By \eqref{eq:2.8}, we have
$$
	\begin{aligned}
	\left|\bigl(g_{1-\alpha,n}*[u-v]_+^2\bigr)(T) - \bigl(g_{1-\alpha}*[u-v]_+^2\bigr)(T)\right|
	&
	\le \|[u-v]_+^2\|_{L^\infty(0,T)}\|g_{1-\alpha,n}-g_{1-\alpha}\|_{L^1(0,T)} \\
	&
	\to 0
	\quad \mbox{as $n\to\infty$},
	\end{aligned}
$$
Moreover, since $r\in L^1(0,T)$ and $f(u) \in L^\infty(0,T)$, by \eqref{eq:2.7}, we get
$$
	\begin{aligned}
	&
	\biggl|\int_0^T\biggl(h_{\alpha,n}*\biggl[r(\cdot)\biggl(f(u(\cdot))-f(v(\cdot))\biggr)\biggr]\biggr)(t)[u(t)-v(t)]_+\,\dee t \\
	& \hspace{5cm}
	- \int_0^Tr(t)\biggl(f(u(t))-f(v(t))\biggr)[u(t)-v(t)]_+\,\dee t\biggr| \\
	&
	\le \left\|\biggl(h_{\alpha,n}*\biggl[r\biggl(f(u)-f(v)\biggr)\biggr]\biggr) - r\biggl(f(u)-f(v)\biggr)\right\|_{L^1(0,T)}\|[u-v]_+\|_{L^\infty(0,T)} \\
	&
	\to 0 \quad \mbox{as $n\to\infty$}.
	\end{aligned}
$$
Therefore, taking the limit $n\to\infty$ in \eqref{eq:3.2}, we obtain
$$
	\frac{1}{2}\bigl(g_{1-\alpha}*[u-v]_+^2\bigr)(T)
	+ \int_0^Tr(t)\biggl(f(u(t))-f(v(t))\biggr)[u(t)-v(t)]_+\,\dee t
	\le 0.
$$
Since $g_{1-\alpha}$, $r$ are nonnegative and $\bigl(f(u)-f(v)\bigr)[u-v]_+ \ge 0$,
by the arbitrariness of $T>0$,
we observe that $[u-v]_+ = 0$ for all $t>0$, that is,
$$
	u(t) \le v(t)
$$
for all $t>0$.
Thus the proof of Lemma~\ref{Lemma:3.1} is complete.
$\Box$
\bigskip

We are now ready to establish two decay estimates for fractional ODEs that will be used in the fast diffusion range.
First, the following result corresponds to the case $0<m<\alpha$.
\begin{lemma}
\label{Lemma:3.2}
	Let $0<\alpha<1$, $\beta > 1-\alpha$, and $\lambda,\mu>0$.
	Assume that $y \in C([0,\infty))$ is the nonnegative solution to the Cauchy problem
	\begin{equation}
	\label{eq:3.3}
		\RL(y-y_0) + \lambda t^{\beta}y = \mu t^{\beta-1}, \quad t>0,
		\qquad y(0)=y_0>0.
	\end{equation}
	Then it holds that
	$$
		y(t) \le Ct^{-1}
	$$
	for sufficiently large $t>0$, where $C>0$ depends only on $\alpha$, $\beta$, $\lambda$, $\mu$, and $y_0$.
\end{lemma}
{\bf Proof.}
We construct a supersolution to problem~\eqref{eq:3.3}, namely, a function $w \in C([0,\infty))$ satisfying
\begin{equation}
\label{eq:3.4}
	\RL(w-y_0)(t) + \lambda t^{\beta}w(t) - \mu t^{\beta-1} \ge 0,
	\quad t>0,
	\qquad w(0)=y_0.
\end{equation}
Define
$$
	t_0^{\alpha+\beta} := \max\left\{\frac{2\Gamma(\alpha+\beta)}{\lambda\Gamma(\beta)}, \frac{1}{\lambda\Gamma(1-\alpha)}\left(2^{1+\alpha}+\frac{2^{2+\alpha}}{1-\alpha}\right)\right\},
$$
and consider the function
\begin{equation}
\label{eq:3.5}
	w(t) =
	\left\{
	\begin{array}{ll}
	y_0 + \dfrac{\mu\Gamma(\beta)}{\Gamma(\alpha+\beta)}t^{\alpha+\beta-1},
	&
	0 \le t < t_0,
	\vspace{5pt}\\
	Kt^{-1},
	&
	t \ge t_0,
	\end{array}
	\right.
\end{equation}
where $K$ are defined by
$$
	K = y_0t_0 + \dfrac{\mu\Gamma(\beta)}{\Gamma(\alpha+\beta)}t_0^{\alpha+\beta}.
$$
Then we easily see that $w \in C([0,\infty)) \cap W_{\loc}^{1,1}([0,\infty))$ and that
$$
	w'(t) =
	\left\{
	\begin{array}{ll}
	\dfrac{\mu\Gamma(\beta)}{\Gamma(\alpha+\beta-1)}t^{\alpha+\beta-2},
	&
	0 \le t < t_0,
	\vspace{5pt}\\
	-Kt^{-2},
	&
	t \ge t_0.
	\end{array}
	\right.
$$
For $t \in [0,t_0)$, it holds that
$$
	\begin{aligned}
	\RL(w-y_0)(t)
	&
	= \bigl(g_{1-\alpha}*w'\bigr)(t) \\
	&
	= \frac{\mu\Gamma(\beta)}{\Gamma(\alpha+\beta-1)\Gamma(1-\alpha)}\int_0^t(t-s)^{-\alpha}s^{\alpha+\beta-2}\,\dee s
	= \mu t^{\beta-1}.
	\end{aligned}
$$
This implies that
$$
	\begin{aligned}
	\RL(w-y_0)(t) + \lambda t^{\beta}w(t) - \mu t^{\beta-1}
	&
	= \mu t^{\beta-1}
	+ \lambda t^{\beta}\left(y_0 + \dfrac{\mu\Gamma(\beta)}{\Gamma(\alpha+\beta)}t^{\alpha+\beta-1}\right)
	-\mu t^{\beta-1} \\
	&
	= \lambda t^{\beta}\left(y_0 + \dfrac{\mu\Gamma(\beta)}{\Gamma(\alpha+\beta)}t^{\alpha+\beta-1}\right)
	\ge 0
	\end{aligned}
$$
for $0\le t<t_0$.
Next, for $t\ge t_0$, we observe that
$$
	\begin{aligned}
	\RL(w-y_0)(t)
	&
	= (g_{1-\alpha}*w')(t) \\
	&
	= \frac{\mu\Gamma(\beta)}{\Gamma(\alpha+\beta-1)}\int_0^{t_0}g_{1-\alpha}(t-s)s^{\alpha+\beta-2}\,\dee s
	- K\int_{t_0}^tg_{1-\alpha}(t-s)s^{-2}\,\dee s \\
	&
	\ge -K\int_{t_0}^tg_{1-\alpha}(t-s)s^{-2}\,\dee s.
	\end{aligned}
$$
If $t_0\le t\le 2t_0$, then we have
$$
	-K\int_{t_0}^tg_{1-\alpha}(t-s)s^{-2}\,\dee s
	\ge -Kt_0^{-2}\frac{(t-t_0)^{1-\alpha}}{\Gamma(2-\alpha)}
	\ge -\frac{K}{\Gamma(2-\alpha)}t_0^{-1-\alpha}
	\ge -\frac{K2^{1+\alpha}}{\Gamma(2-\alpha)}t^{-1-\alpha}.
$$
By the definition of $K$ and $t_0$, it holds that
$$
	K
	\ge \frac{\mu\Gamma(\beta)}{\Gamma(\alpha+\beta)}t_0^{\alpha+\beta}
	\ge \frac{\mu\Gamma(\beta)}{\Gamma(\alpha+\beta)}\frac{2\Gamma(\alpha+\beta)}{\lambda\Gamma(\beta)}
	= \frac{2\mu}{\lambda},
$$
and then we obtain
$$
	\begin{aligned}
	\RL(w-y_0)(t) + \lambda t^{\beta}w(t) - \mu t^{\beta-1}
	&
	\ge -\frac{K2^{1+\alpha}}{\Gamma(2-\alpha)}t^{-1-\alpha} + \lambda K t^{\beta-1} - \frac{\lambda K}{2}t^{\beta-1} \\
	&
	\ge \frac{\lambda K}{2}t^{-\alpha-1}\left(-\frac{2^{2+\alpha}}{\lambda\Gamma(2-\alpha)} + t_0^{\alpha+\beta}\right) \\
	&
	\ge \frac{\lambda K}{2}t^{-\alpha-1}\left(-\frac{2^{2+\alpha}}{\lambda\Gamma(2-\alpha)} + \frac{1}{\lambda\Gamma(1-\alpha)}\left(2^{1+\alpha}+\frac{2^{2+\alpha}}{1-\alpha}\right)\right) \\
	&
	= \frac{K2^{\alpha}}{\Gamma(1-\alpha)}t^{-\alpha-1} \ge 0.
	\end{aligned}
$$
For $t \ge 2t_0$, putting $r = t/t_0 \ge 2$, we deduce that
$$
	\begin{aligned}
	&
	-K\int_{t_0}^tg_{1-\alpha}(t-s)s^{-2}\,\dee s \\
	& \hspace{1cm}
	= -\frac{K}{\Gamma(1-\alpha)}t_0^{-1-\alpha}\int_1^r (r-\tau)^{-\alpha}\tau^{-2}\,\dee\tau \\
	& \hspace{1cm}
	= -\frac{K}{\Gamma(1-\alpha)}t_0^{-1-\alpha}\left(\int_1^{r/2}(r-\tau)^{-\alpha}\tau^{-2}\,\dee\tau + \int_{r/2}^r(r-\tau)^{-\alpha}\tau^{-2}\,\dee\tau\right) \\
	& \hspace{1cm}
	\ge -\frac{K}{\Gamma(1-\alpha)}t_0^{-1-\alpha}\left(2^{\alpha}r^{-\alpha}\int_1^\infty\tau^{-2}\,\dee\tau + 4r^{-2}\int_{r/2}^r(r-\tau)^{-\alpha}\,\dee\tau\right) \\
	& \hspace{1cm}
	= -\frac{K}{\Gamma(1-\alpha)}t_0^{-1-\alpha}\left(2^{\alpha}r^{-\alpha} + \frac{2^{1+\alpha}}{1-\alpha}r^{-1-\alpha}\right) \\
	& \hspace{1cm}
	\ge -\frac{K}{\Gamma(1-\alpha)}t_0^{-1-\alpha}\left(2^{\alpha} + \frac{2^{1+\alpha}}{1-\alpha}\right)r^{\beta-1}
	= -\frac{K}{\Gamma(1-\alpha)}t_0^{-\alpha-\beta}\left(2^{\alpha} + \frac{2^{1+\alpha}}{1-\alpha}\right)t^{\beta-1}.
	\end{aligned}
$$
Here we used the fact that
$$
	r^{-1-\alpha} \le r^{-\alpha} \le r^{\beta-1}
$$
for $r\ge2$ and $\beta-1>-\alpha$.
Therefore, we have
$$
	\begin{aligned}
	&
	\RL(w-y_0)(t) + \lambda t^{\beta}w(t) - \mu t^{\beta-1} \\
	& \hspace{1cm}
	\ge -\frac{K}{\Gamma(1-\alpha)}t_0^{-\alpha-\beta}\left(2^{\alpha} + \frac{2^{1+\alpha}}{1-\alpha}\right)t^{\beta-1}
	+ \frac{\lambda K}{2}t^{\beta-1} \\
	& \hspace{1cm}
	\ge \frac{\lambda K}{2}t_0^{-\alpha-\beta}t^{\beta-1}\left(-\frac{1}{\lambda\Gamma(1-\alpha)}\left(2^{\alpha+1}+\frac{2^{2+\alpha}}{1-\alpha}\right) + t_0^{\alpha+\beta}\right) \\
	& \hspace{1cm}
	\ge 0.
	\end{aligned}
$$
Thus the function $w$ defined by \eqref{eq:3.5} satisfies \eqref{eq:3.4}, so we obtain by Lemma~\ref{Lemma:3.1}
$$
	y(t) \le w(t), \quad t>0.
$$
This implies that
$$
	y(t) \le Ct^{-1}, \quad t\ge t_0.
$$
Hence the proof of Lemma~\ref{Lemma:3.2} is complete.
$\Box$
\bigskip

Finally, we give the following decay estimate for fractional ODEs with logarithmic term.
We will use the following result in the proof of Theorem~\ref{Theorem:1.1} for the critical case $m=\alpha$.
\begin{lemma}
\label{Lemma:3.3}
	Let $0<\alpha<1$, and $\lambda,\mu>0$.
	Assume that $y \in C([0,\infty))$ is the nonnegative solution to
	$$
		\RL(y-y_0) + \lambda t^{1-\alpha}y = \mu t^{-\alpha}\log(1+t), \quad t>0,
		\qquad y(0)=y_0>0.
	$$
	Then it holds that
	$$
		y(t) \le C\frac{\log(1+t)}{t}
	$$
	for sufficiently large $t>0$, where $C>0$ depends only on $\alpha$, $\lambda$, $\mu$, and $y_0$.
\end{lemma}
{\bf Proof.}
As in Lemma~\ref{Lemma:3.1}, we construct a function $z \in C([0,\infty))$ satisfying
\begin{equation}
\label{eq:3.6}
	\RL(z-y_0)(t) + \lambda t^{1-\alpha}z(t) - \mu t^{-\alpha}\log(1+t) \ge 0, \quad t>0.
\end{equation}
Define
$$
	t_0 := \max\left\{e-1, \exp\left(\frac{(2-\alpha)2^{2+\alpha}}{\lambda\Gamma(2-\alpha)}\right)-1, \exp\left(\frac{2^{\alpha+1}}{\lambda\Gamma(1-\alpha)}\right)-1\right\},
$$
and consider the function
\begin{equation}
\label{eq:3.7}
	z(t) =
	\left\{
	\begin{array}{ll}
	y_0 + \mu\Gamma(2-\alpha)t,
	&
	0\le t< t_0,
	\vspace{5pt}\\
	K\dfrac{\log(1+t)}{t},
	&
	t\ge t_0,
	\end{array}
	\right.
\end{equation}
where
$$
	K = \frac{t_0(y_0+\mu\Gamma(2-\alpha)t_0)}{\log(1+t_0)}.
$$
Then we observe that $z \in C([0,\infty)) \cap W_{\loc}^{1,1}([0,\infty))$ and that
$$
	z'(t) =
	\left\{
	\begin{array}{ll}
	\mu\Gamma(2-\alpha),
	&
	0\le t<t_0,
	\vspace{5pt}\\
	K\left(-\dfrac{\log(1+t)}{t^2} + \dfrac{1}{t(1+t)}\right),
	&t\ge t_0.
	\end{array}
	\right.
$$
For $0\le t<t_0$, we have
$$
	\RL(z-y_0)(t)
	= \bigl(g_{1-\alpha}*z'\bigr)(t)
	= \mu t^{1-\alpha}.
$$
Hence, by the elementary inequality
\begin{equation}
\label{eq:3.8}
	\frac{\log(t+1)}{t} \le 1, \quad t\ge0,
\end{equation}
it follows that
$$
	\begin{aligned}
	&
	\RL(z-y_0)(t) + \lambda t^{1-\alpha}z(t) - \mu t^{-\alpha}\log(1+t) \\
	& \hspace{1cm}
	= t^{1-\alpha}\left(\mu + \lambda\left(y_0 + \mu\Gamma(2-\alpha)t\right) - \mu \frac{\log(1+t)}{t}\right) \\
	& \hspace{1cm}
	\ge t^{1-\alpha}\biggl(\mu + \lambda\left(y_0 + \mu\Gamma(2-\alpha)t\right) - \mu\biggr) \\
	& \hspace{1cm}
	= \lambda t^{1-\alpha}\bigl(y_0 + \mu\Gamma(2-\alpha)t\bigr)
	\ge 0.
	\end{aligned}
$$
For $t\ge t_0$, since $g_{1-\alpha} > 0$, by the definition of $z$, we see that
$$
	\begin{aligned}
	\RL(z-y_0)(t)
	&
	= \mu\Gamma(2-\alpha)\int_0^{t_0}g_{1-\alpha}(t-s)\,\dee s 
	+ K\int_{t_0}^tg_{1-\alpha}(t-s)\biggl(-\frac{\log(1+s)}{s^2} + \frac{1}{s(1+s)}\biggr)\,\dee s \\
	&
	\ge -K\int_{t_0}^tg_{1-\alpha}(t-s)\frac{\log(1+s)}{s^2}\,\dee s.
	\end{aligned}
$$
If $t_0\le t\le 2t_0$, using the elementary inequality \eqref{eq:3.8} again, it holds that
$$
	\begin{aligned}
	-K\int_{t_0}^tg_{1-\alpha}(t-s)\frac{\log(1+s)}{s^2}\,\dee s
	&
	\ge -Kt_0^{-1}\int_{t_0}^tg_{1-\alpha}(t-s)\,\dee s \\
	&
	= -Kt_0^{-1}g_{2-\alpha}(t-t_0) \\
	&
	\ge -\frac{K}{\Gamma(2-\alpha)}t_0^{-\alpha}
	\ge -\frac{K2^{\alpha}}{\Gamma(2-\alpha)}t^{-\alpha}.
	\end{aligned}
$$
Note that, by the definition of $t_0$ and the fact that $\exp(s) - 1 \ge s$, $s\in\R$, it follows that
$$
	t_0
	\ge \exp\left(\frac{(2-\alpha)2^{1+\alpha}}{\lambda\Gamma(2-\alpha)}\right) - 1
	\ge \frac{(2-\alpha)2^{1+\alpha}}{\lambda\Gamma(2-\alpha)}
	\ge \frac{2}{\lambda\Gamma(2-\alpha)},
$$
so we get
$$
	K
	\ge \frac{\mu\Gamma(2-\alpha)t_0^2}{\log(1+t_0)}
	\ge \mu\Gamma(2-\alpha)t_0
	\ge \mu\Gamma(2-\alpha)\frac{2}{\lambda\Gamma(2-\alpha)}
	= \frac{2\mu}{\lambda}.
$$
Thus we have for $t_0 \le t \le 2t_0$
$$
	\begin{aligned}
	&
	\RL(z-y_0)(t) + \lambda t^{1-\alpha}z(t) - \mu t^{-\alpha}\log(1+t) \\
	& \hspace{1cm}
	\ge -\frac{K2^{\alpha}}{\Gamma(2-\alpha)}t^{-\alpha} + \frac{\lambda K}{2}t^{-\alpha}\log(1+t) \\
	& \hspace{1cm}
	\ge \frac{\lambda K}{2}t^{-\alpha}\left(-\frac{2^{\alpha+1}}{\lambda\Gamma(2-\alpha)} + \log(1+t_0)\right) \ge 0.
	\end{aligned}
$$
For $t \ge 2t_0$, putting $r = t/t_0 \ge 2$, we deduce that
$$
	\begin{aligned}
	&
	-K\int_{t_0}^tg_{1-\alpha}(t-s)\frac{\log(1+s)}{s^2}\,\dee s \\
	&
	= -\frac{K}{\Gamma(1-\alpha)}t_0^{-\alpha-1}\int_1^r(r-\tau)^{-\alpha}\frac{\log(1+t_0\tau)}{\tau^2}\,\dee\tau \\
	&
	= -\frac{K}{\Gamma(1-\alpha)}t_0^{-\alpha-1}\left(\int_1^{r/2}(r-\tau)^{-\alpha}\frac{\log(1+t_0\tau)}{\tau^2}\,\dee\tau + \int_{r/2}^r(r-\tau)^{-\alpha}\frac{\log(1+t_0\tau)}{\tau^2}\,\dee\tau\right).
	\end{aligned}
$$
For $\tau \in (1,r/2)$, since $\tau\ge1$ and $1+t_0 \ge e$, it holds that
$$
	\log(1+t_0\tau) \le \log\bigl((1+t_0)\tau\bigr)
$$
and that $\log\bigl((1+t_0)\tau\bigr) \ge 1$.
Then we see that
\begin{equation}
\label{eq:3.9}
	\begin{aligned}
	\int_1^{r/2}(r-\tau)^{-\alpha}\frac{\log(1+t_0\tau)}{\tau^2}\,\dee\tau
	&
	 \le 2^{\alpha}r^{-\alpha}\int_1^{r/2}\frac{\log(1+t_0\tau)}{\tau^2}\,\dee\tau \\
	&
	\le 2^{\alpha}r^{-\alpha}\left(\int_1^{r/2}\frac{\log\bigl((1+t_0)\tau\bigr)-1}{\tau^2}\,\dee\tau + \int_1^{r/2}\tau^{-2}\,\dee\tau\right). 
	\end{aligned}
\end{equation}
Note that
$$
	0 \le \frac{\log\bigl((1+t_0)\tau\bigr)-1}{\tau^2} = -\dtau\left(\frac{\log\bigl((1+t_0)\tau\bigr)}{\tau}\right),
$$
we have
$$
	\int_1^{r/2}\frac{\log\bigl((1+t_0)\tau\bigr)-1}{\tau^2}\,\dee\tau
	\le -\int_1^\infty\dtau\left(\frac{\log\bigl((1+t_0)\tau\bigr)}{\tau}\right)\,\dee\tau
	= \log(1+t_0).
$$
Moreover, we also observe that
$$
	\int_1^{r/2}\tau^{-2}\,\dee\tau
	\le \int_1^\infty\tau^{-2}\,\dee\tau
	= 1.
$$
These together with \eqref{eq:3.9} imply that
$$
	\int_1^{r/2}(r-\tau)^{-\alpha}\frac{\log(1+t_0\tau)}{\tau^2}\,\dee\tau
	\le 2^{\alpha}r^{-\alpha}\log(1+t_0) + 2^{\alpha}r^{-\alpha}.
$$
Since $\log\bigl((1+t_0)r\bigr)/r$ is monotone decreasing for $1+t_0 \ge e$ and $r\ge1$, we observe that
$$
	\frac{\log\bigl((1+t_0)r\bigr)}{r} \le \log(1+t_0), \quad r\ge1.
$$
Thus, for $\tau \in (r/2,r)$, we obtain
$$
	\begin{aligned}
	\int_{r/2}^r(r-\tau)^{-\alpha}\frac{\log(1+t_0\tau)}{\tau^2}\,\dee\tau
	&
	\le 4r^{-2}\log\bigl((1+t_0)r\bigr)\int_{r/2}^r(r-\tau)^{-\alpha}\,\dee\tau \\
	&
	\le \frac{2^{1+\alpha}}{1-\alpha}r^{-\alpha}\log(1+t_0).
	\end{aligned}
$$
Hence, since $r=t/t_0$ and $t_0 \ge 1$, we have
$$
	\begin{aligned}
	\RL(z-y_0)(t)
	&
	\ge -\frac{K}{\Gamma(1-\alpha)}t_0^{-\alpha-1}\int_1^r(r-\tau)^{-\alpha}\frac{\log(1+t_0\tau)}{\tau^2}\,\dee\tau \\
	&
	\ge -\frac{K}{\Gamma(1-\alpha)}t_0^{-\alpha-1}\left(2^{\alpha}r^{-\alpha}\log(1+t_0) + 2^{\alpha}r^{-\alpha} + \frac{2^{1+\alpha}}{1-\alpha}r^{-\alpha}\log(1+t_0)\right) \\
	&
	= -\frac{K}{\Gamma(1-\alpha)}t^{-\alpha}\left(2^{\alpha}\frac{\log(1+t_0)}{t_0} + \frac{2^{\alpha}}{t_0} + \frac{2^{1+\alpha}}{1-\alpha}\frac{\log(1+t_0)}{t_0}\right) \\
	&
	\ge -\frac{K}{\Gamma(1-\alpha)}t^{-\alpha}\left(2^{\alpha} + 2^{\alpha} + \frac{2^{\alpha+1}}{1-\alpha}\right) \\
	&
	= -\frac{K(2-\alpha)2^{1+\alpha}}{\Gamma(2-\alpha)}t^{-\alpha}.
	\end{aligned}
$$
Furthermore, we get
$$
	\begin{aligned}
	\RL(z-y_0)(t) + \lambda t^{1-\alpha}z(t) - \mu t^{-\alpha}\log(1+t)
	&
	\ge -\frac{K(2-\alpha)}{\Gamma(2-\alpha)}2^{1+\alpha}t^{-\alpha} + \frac{\lambda K}{2}t^{-\alpha}\log(1+t) \\
	&
	\ge \frac{\lambda K}{2}t^{-\alpha}\left(-\frac{(2-\alpha)2^{2+\alpha}}{\lambda\Gamma(2-\alpha)}+\log(1+t_0)\right) \\
	&
	\ge 0
	\end{aligned}
$$
for $t\ge2t_0$.
Hence the function $z$ defined by \eqref{eq:3.7} satisfies \eqref{eq:3.6}.
Then, by Lemma~\ref{Lemma:3.1}, it holds that
$$
	y(t) \le z(t), \quad t>0.
$$
Thus we conclude that
$$
	y(t) \le C\frac{\log(1+t)}{t}
$$
for $t\ge t_0$, so the proof of Lemma~\ref{Lemma:3.3} is complete.
$\Box$

\section{Proof of Theorem~\ref{Theorem:1.1}}
\label{Section:4}

In this section, we prove Theorem~\ref{Theorem:1.1}.
\subsection{A basic entropy inequality}
\label{Subsection:4.1}
We first prove the following entropy estimate.
\begin{lemma}
\label{Lemma:4.1}
	Assume the same conditions as in Theorem~$\ref{Theorem:1.1}$.
	Let $T>0$ be arbitrary, and let $v(x,t) := t^{\frac{\alpha}{m}}u(x,t)$ and the functional $\EE[v|S]$ be defined by
	$$
		\EE[v|S](t) := \int_{\Omega}\left\{\frac{v^{m+1}(x,t)}{m+1} - \frac{S^{m+1}(x)}{m+1} - S^m(x)\bigl(v(x,t)-S(x)\bigr)\right\}\,\dee x
	$$
	with
	$$
		\EE_0 := \EE[v|S](0) = \int_{\Omega}\left\{-\frac{S^{m+1}(x)}{m+1} + S^{m+1}(x)\right\}\,\dee x = \frac{m}{m+1}\int_{\Omega}S^{m+1}(x)\,\dee x.
	$$
	Then the following assertions hold.
	\begin{itemize}
	\item[$(a)$]
	$\EE := \EE[v|S]$
	satisfies
	\begin{equation}
	\label{eq:4.1}
		\dt\bigl(g_{1-\alpha,n}*[\EE-\EE_0]\bigr)(t)
		+ t^{\alpha\left(\frac{1}{m}-1\right)}\|\nabla(v^m(t)-S^m)\|_2^2
		\le F(t)
		+ I_n(t) + J_n(t)
	\end{equation}
	for a.e. $0<t<T$, where
	$$
		F(t) := C\int_0^t(t-s)^{-\alpha-1}\left(t^{\frac{\alpha}{m}}-s^{\frac{\alpha}{m}}\right)\|u(s)\|_{m+1}\,\dee s,
	$$
	and $C>0$ depends only on $m$, $\alpha$, $N$, $\Omega$, and $\|u_0\|_{p_0}$.
	Moreover, $I_n, J_n \to 0$ strongly in $L^1(0,T)$ as $n\to\infty$.
	
	\item[$(b)$]
	$F$ satisfies the following estimates
	\begin{equation}
	\label{eq:4.2}
		F(t)
		\le
		\left\{
		\begin{array}{lll}
		\nu g_{1-\alpha}(t)
		&
		\mbox{if} \quad m>\alpha,
		\vspace{5pt}\\
		\mu t^{-\alpha}\log(1+t) + \nu g_{1-\alpha}(t)
		&
		\mbox{if} \quad m=\alpha,
		\vspace{5pt}\\
		\mu t^{\beta-1} + \nu g_{1-\alpha}(t)
		&
		\mbox{if} \quad 0<m<\alpha,
		\end{array}
		\right.
		\quad 0<t<T,
	\end{equation}
	where $\beta = \alpha(1/m-1)$, and $\mu$, $\nu$ are positive constants depending only on
	$m$, $\alpha$, $N$, $\Omega$, and $\|u_0\|_{L^{p_0}(\Omega)}$.
	\end{itemize}
\end{lemma}
\bigskip

The estimates \eqref{eq:4.2} imply that the fast diffusion range $0<m<1$ must be divided into three cases: $\alpha<m<1$, $m=\alpha$, and $0<m<\alpha$.
In particular, the critical case $m=\alpha$ produces the logarithmic term.
In order to prove Lemma~\ref{Lemma:4.1},
we recall the sharp decay estimate for $H^{-1}$-solutions to problem~\eqref{P}.
\begin{proposition}
\label{Proposition:4.1}
\cite{BGI}*{Theorem~1.7}
	Assume the same conditions as in Theorem~$\ref{Theorem:1.1}$.
	Then, for every $1\le p\le\infty$,
	\begin{equation}
	\label{eq:4.3}
		\frac{c_0}{1+t^{\frac{\alpha}{m}}} \le \|u(t)\|_p \le \frac{c_1}{t^{\frac{\alpha}{m}}}, \quad t>0,
	\end{equation}
	where $c_0,c_1>0$ depend only on $N$, $\alpha$, $p$, $m$, $\Omega$, and $\|u_0\|_{p_0}$.
\end{proposition}
\bigskip

We are now ready to prove Lemma~\ref{Lemma:4.1}.

\noindent
{\bf Proof of assertion $(a)$ of Lemma~\ref{Lemma:4.1}.}
Let $h_{\alpha,n}$ be the nonnegative kernel defined by \eqref{eq:2.6}.
Convolving \eqref{eq:1.1} with the kernel $h_{\alpha,n}$ and using the identity
$$
	\biggl(h_{\alpha,n}*\biggl[\dt\bigl(g_{1-\alpha}*[u-u_0]\bigr)\biggr]\biggr)(t)
	= \dt\bigl(g_{1-\alpha,n}*[u-u_0]\bigr)(t),
	\quad \mbox{a.e. $0<t<T$},
$$
we obtain
\begin{equation}
\label{eq:4.4}
	\begin{aligned}
	&
	\int_{\Omega}\dt\bigl(g_{1-\alpha,n}*u(x,\cdot)\bigr)(t)\varphi(x)\,\dee x \\
	& \hspace{1cm}
	+ \int_{\Omega}\biggl(h_{\alpha,n}*\nabla u^m(x,\cdot)\biggr)(t)\cdot\nabla\varphi(x)\,\dee x
	- g_{1-\alpha,n}(t)\int_{\Omega}u_0(x)\varphi(x)\,\dee x
	= 0
	\end{aligned}
\end{equation}
for a.e. $t\in(0,T)$.
Since $g_{1-\alpha,n} \in W_{\loc}^{1,1}([0,\infty))$ and $u \in L^\infty(0,\infty; H^{-1}(\Omega))$,
by Fubini's theorem, we get
$$
	\begin{aligned}
	&
	\int_{\Omega}\dt\bigl(g_{1-\alpha,n}*u(x,\cdot)\bigr)(t)\varphi(x)\,\dee x \\
	& \hspace{1cm}
	= \int_{\Omega}\biggl(\bigl(g_{1-\alpha,n}'*u(x,\cdot)\bigr)(t) + g_{1-\alpha,n}(0)u(x,t)\biggr)\varphi(x)\,\dee x \\
	& \hspace{1cm}
	= \biggl(g_{1-\alpha,n}'*\biggl[\int_{\Omega}u(x,\cdot)\varphi(x)\,\dee x\biggr]\biggr)(t) + g_{1-\alpha,n}(0)\int_{\Omega}u(x,t)\varphi(x)\,\dee x \\
	& \hspace{1cm}
	= \dt\biggl(g_{1-\alpha,n}*\biggl[\int_{\Omega}u(x,\cdot)\varphi(x)\,\dee x\biggr]\biggr)(t), \quad \mbox{a.e. $0<t<T$}.
	\end{aligned}
$$
Multiplying \eqref{eq:4.4} by $t^{\frac{\alpha}{m}}$ and applying this, we deduce that
\begin{equation}
\label{eq:4.5}
	\begin{aligned}
	&
	t^{\frac{\alpha}{m}}\dt\biggl(g_{1-\alpha,n}*\biggl[\int_{\Omega}u(x,\cdot)\varphi(x)\,\dee x\biggr]\biggr)(t) \\
	& \hspace{1cm}
	+ t^{\frac{\alpha}{m}}\int_{\Omega}\biggl(h_{\alpha,n}*\nabla u^m(x,\cdot)\biggr)(t)\cdot\nabla\varphi(x)\,\dee x
	- t^{\frac{\alpha}{m}}g_{1-\alpha,n}(t)\int_{\Omega}u_0(x)\varphi(x)\,\dee x
	= 0
	\end{aligned}
\end{equation}
for a.e. $t\in(0,T)$.
Since
$$
	\int_{\Omega}u(x,\cdot)\varphi(x)\,\dee x \in L^1(0,T),
$$
for the first term on the left-hand side,
by using \eqref{eq:2.10}, we get
$$
	\begin{aligned}
	&
	t^{\frac{\alpha}{m}}\dt\biggl(g_{1-\alpha,n}*\biggl[\int_{\Omega}u(x,\cdot)\varphi(x)\,\dee x\biggr]\biggr)(t) \\
	&
	= \dt\biggl(g_{1-\alpha,n}*\biggl[t^{\frac{\alpha}{m}}\int_{\Omega}u(x,\cdot)\varphi(x)\,\dee x\biggr]\biggr)(t) \\
	& \hspace{4cm}
	+ \int_0^tg_{1-\alpha,n}'(t-s)\left(t^{\frac{\alpha}{m}} - s^{\frac{\alpha}{m}}\right)\left\{\int_{\Omega}u(x,s)\varphi(x)\,\dee x\right\}\,\dee s \\
	&
	= \dt\biggl(g_{1-\alpha,n}*\biggl[\int_{\Omega}v(x,\cdot)\varphi(x)\,\dee x\biggr]\biggr)(t) \\
	& \hspace{4cm}
	+ \int_0^tg_{1-\alpha,n}'(t-s)\left(t^{\frac{\alpha}{m}} - s^{\frac{\alpha}{m}}\right)\left\{\int_{\Omega}u(x,s)\varphi(x)\,\dee x\right\}\,\dee s
	\end{aligned}
$$
for a.e. $0<t<T$.
Using Fubini's theorem again, we obtain
$$
	\dt\biggl(g_{1-\alpha,n}*\biggl[\int_{\Omega}v(x,\cdot)\varphi(x)\,\dee x\biggr]\biggr)(t)
	= \int_{\Omega}\dt\bigl(g_{1-\alpha,n}*v(x,\cdot)\bigr)(t)\varphi(x)\,\dee x
$$
for a.e. $0<t<T$.
These together with \eqref{eq:4.5} imply
$$
	\begin{aligned}
	&
	\int_{\Omega}\dt\bigl(g_{1-\alpha,n}*v(x,\cdot)\bigr)(t)\varphi(x)\,\dee x \\
	& \hspace{1cm}
	+ t^{\frac{\alpha}{m}}\int_{\Omega}\biggl(h_{\alpha,n}*\nabla u^m(x,\cdot)\biggr)(t)\cdot\nabla\varphi(x)\,\dee x
	- t^{\frac{\alpha}{m}}g_{1-\alpha,n}(t)\int_{\Omega}u_0(x)\varphi(x)\,\dee x \\
	&
	= -\int_0^tg_{1-\alpha,n}'(t-s)\left(t^{\frac{\alpha}{m}} - s^{\frac{\alpha}{m}}\right)\left\{\int_{\Omega}u(x,s)\varphi(x)\,\dee x\right\}\,\dee s,
	\quad \mbox{a.e. $0<t<T$}.
	\end{aligned}
$$
We choose $\varphi = v^m(t) - S^m \in H_0^1(\Omega)$ as a test function to get
\begin{equation}
\label{eq:4.6}
	\begin{aligned}
	&
	\int_{\Omega}\dt\bigl(g_{1-\alpha,n}*v(x,\cdot)\bigr)(t)(v^m(x,t) - S^m(x))\,\dee x \\
	& \hspace{1cm}
	+ t^{\frac{\alpha}{m}}\int_{\Omega}\biggl(h_{\alpha,n}*\nabla u^m(x,\cdot)\biggr)(t)\cdot\nabla(v^m(x,t) - S^m(x))\,\dee x \\
	& \hspace{3cm}
	- t^{\frac{\alpha}{m}}g_{1-\alpha,n}(t)\int_{\Omega}u_0(x)(v^m(x,t) - S^m(x))\,\dee x \\
	&
	= -\int_0^tg_{1-\alpha,n}'(t-s)\left(t^{\frac{\alpha}{m}} - s^{\frac{\alpha}{m}}\right)\left\{\int_{\Omega}u(x,s)(v^m(x,t) - S^m(x))\,\dee x\right\}\,\dee s
	\end{aligned}
\end{equation}
for a.e. $t\in(0,T)$.
We now write
$$
	\begin{aligned}
	&
	t^{\frac{\alpha}{m}}\int_{\Omega}\biggl(h_{\alpha,n}*\nabla u^m(x,\cdot)\biggr)(t)\cdot\nabla(v^m(x,t) - S^m(x))\,\dee x \\
	&
	= t^{\frac{\alpha}{m}}\int_{\Omega}\nabla u^m(x,t)\cdot\nabla(v^m(x,t)-S^m(x))\,\dee x \\
	& \hspace{2cm}
	- t^{\frac{\alpha}{m}}\int_{\Omega}\biggl\{\nabla u^m(x,t)-\biggl(h_{\alpha,n}*\nabla u^m(x,\cdot)\biggr)(t)\biggr\}\cdot\nabla(v^m(x,t)-S^m(x))\,\dee x \\
	&
	=: t^{\frac{\alpha}{m}}\int_{\Omega}\nabla u^m(x,t)\cdot\nabla(v^m(x,t)-S^m(x))\,\dee x - I_n(t),
	\end{aligned}
$$
and
$$
	\begin{aligned}
	&
	t^{\frac{\alpha}{m}}g_{1-\alpha,n}(t)\int_{\Omega}u_0(x)(v^m(x,t) - S^m(x))\,\dee x \\
	&
	= t^{\frac{\alpha}{m}}g_{1-\alpha}(t)\int_{\Omega}u_0(x)(v^m(x,t) - S^m(x))\,\dee x \\
	& \hspace{2cm}
	- t^{\frac{\alpha}{m}}\bigl(g_{1-\alpha}(t) - g_{1-\alpha,n}(t)\bigr)\int_{\Omega}u_0(x)(v^m(x,t)-S^m(x))\,\dee x \\
	&
	=: \frac{t^{\alpha\left(\frac{1}{m}-1\right)}}{\Gamma(1-\alpha)}\int_{\Omega}u_0(x)(v^m(x,t) - S^m(x))\,\dee x - J_n(t).
	\end{aligned}
$$
Then it holds that
$$
	I_n, J_n \to 0 \quad \mbox{strongly in $L^1(0,T)$ as $n\to\infty$}.
$$
Indeed, by H\"{o}lder's inequality and \eqref{eq:2.7}, we obtain
$$
	\begin{aligned}
	\int_0^T|I_n(t)|\,\dee t
	&
	\le T^{\frac{\alpha}{m}}\|\nabla u^m - (h_{\alpha,n}*\nabla u^m)\|_{L^2(0,T; L^2(\Omega))}\|\nabla(v^m-S^m)\|_{L^2(0,T; L^2(\Omega))} \\
	&
	\to 0 \quad \mbox{as $n\to\infty$}.
	\end{aligned}
$$
Note that, by \eqref{eq:4.3} in Proposition~\ref{Proposition:4.1} and assertion (ii) of Remark~\ref{Remark:1.1},
it holds that
$$
	\|v(t)\|_{\infty} = t^{\frac{\alpha}{m}}\|u(t)\|_{\infty} \le C,
	\qquad \|S\|_{\infty} \le C\|u_0\|_{p_0}^{\frac{1}{m}}.
$$
Since $\Omega$ is bounded and $u_0 \in L^1(\Omega)$,
these estimates together with \eqref{eq:2.8} imply that
$$
	\int_0^T|J_n(t)|\,\dee t
	\le CT^{\frac{\alpha}{m}}\|g_{1-\alpha,n}-g_{1-\alpha}\|_{L^1(0,T)}
	\to 0 \quad \mbox{as $n\to\infty$}.
$$
Moreover, by \eqref{eq:1.6} with $\varphi = v^m(t)-S^m$, we have
$$
	\frac{t^{\alpha\left(\frac{1}{m}-1\right)}}{\Gamma(1-\alpha)}\int_{\Omega}u_0(x)(v^m(x,t) - S^m(x))\,\dee x
	= t^{\alpha\left(\frac{1}{m}-1\right)}\int_{\Omega}\nabla S^m(x)\cdot\nabla(v^m(x,t)-S^m(x))\,\dee x
$$
for a.e. $t>0$.
Hence, since $u^m = t^{-\alpha}v^m$, we can rewrite \eqref{eq:4.6} as
\begin{equation}
\label{eq:4.7}
	\begin{aligned}
	&
	\int_{\Omega}\dt\bigl(g_{1-\alpha,n}*v(x,\cdot)\bigr)(t)(v^m(x,t)-S^m(x))\,\dee x 
	+ t^{\alpha\left(\frac{1}{m}-1\right)}\int_{\Omega}|\nabla(v^m(x,t)-S^m(x))|^2\,\dee x \\
	&
	= -\int_0^tg_{1-\alpha,n}'(t-s)\left(t^{\frac{\alpha}{m}} - s^{\frac{\alpha}{m}}\right)\left\{\int_{\Omega}u(x,s)(v^m(x,t) - S^m(x))\,\dee x\right\}\,\dee s
	+ I_n(t) + J_n(t)
	\end{aligned}
\end{equation}
for a.e. $t\in(0,T)$.
For the first term in the left-hand side, using \eqref{eq:2.9} for $H(y) = \frac{1}{m+1}y^{m+1}$, we deduce that
$$
	\frac{1}{m+1}\dt\bigl(g_{1-\alpha,n}*v^{m+1}(x,\cdot)\bigr)(t)
	\le v^m(x,t)\dt\bigl(g_{1-\alpha,n}*v(x,\cdot)\bigr)(t), \quad \mbox{a.e. $x\in\Omega$, $t\in(0,T)$}.
$$
Since
$$
	\begin{aligned}
	\EE(t) - \EE_0
	&
	= \int_{\Omega}\left\{\frac{v^{m+1}(x,t)}{m+1} - \frac{S^{m+1}(x)}{m+1} - S^m(x)\bigl(v(x,t)-S(x)\bigr)\right\}\,\dee x \\
	& \hspace{5cm}
	- \int_{\Omega}\left\{-\frac{S^{m+1}(x)}{m+1} + S^{m+1}(x)\right\}\,\dee x \\
	&
	= \int_{\Omega}\left\{\frac{v^{m+1}(x,t)}{m+1} - S^m(x)v(x,t)\right\}\,\dee x,
	\end{aligned}
$$
we have
$$
	\begin{aligned}
	&
	\dt\bigl(g_{1-\alpha,n}*[\EE-\EE_0]\bigr)(t) \\
	&
	= \int_{\Omega}\left\{\frac{1}{m+1}\dt\bigl(g_{1-\alpha,n}*v^{m+1}(x,\cdot)\bigr)(t) - S^m(x)\dt\bigl(g_{1-\alpha,n}*v(x,\cdot)\bigr)(t)\right\}\,\dee x \\
	&
	\le \int_{\Omega}\left\{v^m(x,t)\dt\bigl(g_{1-\alpha,n}*v(x,\cdot)\bigr)(t) - S^m(x)\dt\bigl(g_{1-\alpha,n}*v(x,\cdot)\bigr)(t)\right\}\,\dee x \\
	&
	= \int_{\Omega}\dt\bigl(g_{1-\alpha,n}*v(x,\cdot)\bigr)(t)(v^m(x,t)-S^m(x))\,\dee x,
	\quad \mbox{a.e. $0<t<T$}.
	\end{aligned}
$$
This together with \eqref{eq:4.7} implies
\begin{equation}
\label{eq:4.8}
	\begin{aligned}
	&
	\dt\bigl(g_{1-\alpha,n}*[\EE-\EE_0]\bigr)(t)
	+ t^{\alpha\left(\frac{1}{m}-1\right)}\|\nabla(v^m(t)-S^m)\|_2^2 \\
	& \hspace{1cm}
	\le -\int_0^tg_{1-\alpha,n}'(t-s)\left(t^{\frac{\alpha}{m}}-s^{\frac{\alpha}{m}}\right)\left\{\int_{\Omega}u(x,s)(v^m(x,t)-S^m(x))\,\dee x\right\}\,\dee s
	\end{aligned}
\end{equation}
for a.e. $t\in(0,T)$.
We next estimate the first term on the right-hand side of \eqref{eq:4.8}.
From \eqref{eq:2.2} and \eqref{eq:2.3},
since $-g_{1-\alpha,n}' \ge 0$, we obtain the following uniform estimate with respect to $n$
$$
	\begin{aligned}
	-g_{1-\alpha,n}'(t)
	&
	=-ns_{\alpha,n}'(t)
	= n^2t^{\alpha-1}E_{\alpha,\alpha}(-nt^{\alpha})
	\le \frac{1}{\Gamma(\alpha)}\frac{n^2t^{\alpha-1}}{(1+c_{\alpha,2}nt^{\alpha})^2} \\
	&
	= t^{-\alpha-1}\frac{(nt^{\alpha})^2}{\Gamma(\alpha)(1+c_{\alpha,2}nt^{\alpha})^2}
	\le C_{\alpha}t^{-\alpha-1}, \quad 0<t<T,
	\end{aligned}
$$
where $C_{\alpha}$ depends only on $\alpha$.
Since $\Omega$ is bounded, by the boundedness of $v$ and $S$, we also deduce that
$$
	\begin{aligned}
	&
	\left|\int_0^t\left[-g_{1-\alpha,n}'(t-s)\right]\bigl(t^{\frac{\alpha}{m}}-s^{\frac{\alpha}{m}}\bigr)\left\{\int_{\Omega}(v^m(x,t)-S^m(x))u(x,s)\,\dee x\right\}\,\dee s\right| \\
	&
	\le C\int_0^t\left[-g_{1-\alpha,n}'(t-s)\right]\bigl(t^{\frac{\alpha}{m}}-s^{\frac{\alpha}{m}}\bigr)\|u(s)\|_{m+1}\,\dee s \\
	&
	\le C\int_0^t(t-s)^{-\alpha-1}\left(t^{\frac{\alpha}{m}}-s^{\frac{\alpha}{m}}\right)\|u(s)\|_{m+1}\,\dee s,
	\quad \mbox{a.e. $0<t<T$}.
	\end{aligned}
$$
This together with \eqref{eq:4.8} yields \eqref{eq:4.1}.
Therefore, assertion $(a)$ of Lemma~\ref{Lemma:4.1} follows.
$\Box$
\bigskip

\noindent
{\bf Proof of assertion $(b)$ of Lemma~\ref{Lemma:4.1}}.
Combining \eqref{eq:4.3} in Proposition~\ref{Proposition:4.1} and \eqref{eq:1.2},
we have
$$
	\|u(t)\|_{m+1} \lesssim (1+t)^{-\frac{\alpha}{m}}, \quad t>0,
$$
where $a\lesssim b$ means that $a \le Cb$ for some constant $C$ independent of $a$ and $b$.
If $m>\alpha$, then $\alpha/m<1$, and \eqref{eq:4.3} in Proposition~\ref{Proposition:4.1} also gives
$\|u(t)\|_{m+1}\lesssim t^{-\frac{\alpha}{m}}$ and $t^{-\frac{\alpha}{m}} \in L_{\loc}^1([0,\infty))$.
Then, since $t^{-\alpha} = \Gamma(1-\alpha)g_{1-\alpha}(t)$, we obtain
$$
	\begin{aligned}
	F(t)
	&
	= C\int_0^{t/2}(t-s)^{-\alpha-1}\left(t^{\frac{\alpha}{m}} - s^{\frac{\alpha}{m}}\right)\|u(s)\|_{m+1}\,\dee s \\
	& \hspace{3cm}
	+ C\int_{t/2}^t(t-s)^{-\alpha-1}\left(t^{\frac{\alpha}{m}} - s^{\frac{\alpha}{m}}\right)\|u(s)\|_{m+1}\,\dee s \\
	&
	\lesssim t^{-\alpha-1+\frac{\alpha}{m}}\int_0^{t/2}s^{-\frac{\alpha}{m}}\,\dee s
	+ t^{\frac{\alpha}{m}-1}\int_{t/2}^t(t-s)^{-\alpha}s^{-\frac{\alpha}{m}}\,\dee s
	\lesssim  g_{1-\alpha}(t).
	\end{aligned}
$$
On the other hand, when $0<m\le\alpha$, we use the global estimate with $(1+s)^{-\alpha/m}$.
If $m=\alpha$, then it follows that
$$
	\begin{aligned}
	F(t)
	&
	\lesssim t^{-\alpha}\int_0^{t/2}(1+s)^{-1}\,\dee s
	+ t^{-1}\int_{t/2}^t(t-s)^{-\alpha}\,\dee s\\
	&
	\lesssim t^{-\alpha}\log(1+t) + g_{1-\alpha}(t).
	\end{aligned}
$$
If $0<m<\alpha$, setting
$$
	\beta:=\alpha\left(\frac1m-1\right)>1-\alpha,
$$
we similarly obtain
$$
	\begin{aligned}
	F(t)
	&
	\lesssim t^{-\alpha-1+\frac{\alpha}{m}}\int_0^{t/2}(1+s)^{-\frac{\alpha}{m}}\,\dee s
	+ t^{\frac{\alpha}{m}-1}\int_{t/2}^t(t-s)^{-\alpha}(1+s)^{-\frac{\alpha}{m}}\,\dee s\\
	&
	\lesssim t^{\beta-1} + g_{1-\alpha}(t).
	\end{aligned}
$$
Thus, assertion $(b)$ of Lemma~\ref{Lemma:4.1} follows,
and the proof of Lemma~\ref{Lemma:4.1} is complete.
$\Box$

\subsection{The porous medium and linear cases $m\ge1$}
\label{Subsection:4.2}
We first prove Theorem~\ref{Theorem:1.1} for the case $m\ge1$.
We first establish the following estimate
\begin{equation}
\label{eq:4.9}
	\lambda\EE^m(t)
	\le \|\nabla(v^m(t)-S^m)\|_2^2,
	\quad \mbox{a.e. $0<t<T$},
\end{equation}
where $\lambda > 0$ depends only on $N$, $\alpha$, $m$, $\Omega$, and $\|u_0\|_{p_0}$.
By Taylor's theorem, there exists a $c\in(a,b)$ such that
$$
	\frac{a^{m+1}}{m+1} - \frac{b^{m+1}}{m+1} - b^m(a-b) = \frac{m}{2}c^{m-1}(a-b)^2.
$$
Hence there exists a nonnegative function
$\eta \in L^\infty(\Omega\times(0,\infty))$ satisfying
$$
	0 \le \eta(x,t) \le M := \max\left\{\esssup_{t>0}\|v(t)\|_\infty, \|S\|_\infty\right\} < \infty
$$
and
$$
	\frac{v^{m+1}(x,t)}{m+1} - \frac{S^{m+1}(x)}{m+1} - S^m(x)(v(x,t)-S(x))
	= \frac{m}{2}\eta^{m-1}(x,t)|v(x,t)-S(x)|^2.
$$
This together with $m\ge1$ yields
$$
	\begin{aligned}
	\EE(t)
	&
	= \int_{\Omega}\left\{\frac{v^{m+1}(x,t)}{m+1} - \frac{S^{m+1}(x)}{m+1} - S^m(x)(v(x,t)-S(x))\right\}\,\dee x \\
	&
	= \frac{m}{2}\int_{\Omega}\eta^{m-1}(x,t)|v(x,t)-S(x)|^2\,\dee x \\
	&
	\le \frac{mM^{m-1}}{2}\int_{\Omega}|v(x,t)-S(x)|^2\,\dee x,
	\quad \mbox{a.e. $0<t<T$}.
	\end{aligned}
$$
Furthermore, by using the elementary inequality
$$
	|a-b| \le |a^m-b^m|^{\frac{1}{m}},
	\quad a,b \ge 0,
	\quad m\ge1,
$$
H\"{o}lder's inequality and the Poincar\'{e} inequality
$$
	\begin{aligned}
	\int_{\Omega}|v(x,t)-S(x)|^2\,\dee x
	&
	\le \int_{\Omega}|v^m(x,t)-S^m(x)|^{\frac{2}{m}}\,\dee x \\
	&
	\le \left(\int_{\Omega}\,\dee x\right)^{1-\frac{1}{m}}\left(\int_{\Omega}|v^m(x,t)-S^m(x)|^2\,\dee x\right)^{\frac{1}{m}} \\
	&
	= C\|v^m(t)-S^m\|_2^{\frac{2}{m}} \\
	&
	\le C\|\nabla(v^m(t)-S^m)\|_2^{\frac{2}{m}},
	\quad \mbox{a.e. $0<t<T$}.
	\end{aligned}
$$
Thus we get
$$
	\EE(t) \le C\|\nabla(v^m(t)-S^m)\|_2^{\frac{2}{m}}, \quad \mbox{a.e. $0<t<T$},
$$
that is, \eqref{eq:4.9}.
Combining \eqref{eq:4.1}, \eqref{eq:4.9} and the first estimate in \eqref{eq:4.2},
we have 
$$
	\dt\bigl(g_{1-\alpha,n}*[\EE - \EE_0]\bigr)(t)
	+ \lambda t^{-\alpha\left(1-\frac{1}{m}\right)}\EE^m(t)
	\le \nu g_{1-\alpha}(t) + I_n(t) + J_n(t),
	\quad \mbox{a.e. $0<t<T$}.
$$
Moreover, since
$$
	g_{1-\alpha,n}(t) = \dt\bigl(g_{1-\alpha,n}*1\bigr)(t),
$$
we can write
$$
	\begin{aligned}
	\nu g_{1-\alpha}(t)
	&
	= \nu g_{1-\alpha,n}(t) + \nu\bigl(g_{1-\alpha}(t)-g_{1-\alpha,n}(t)\bigr) \\
	&
	= \nu\dt\bigl(g_{1-\alpha,n}*1\bigr)(t) + \nu\bigl(g_{1-\alpha}(t)-g_{1-\alpha,n}(t)\bigr) \\
	&
	=: \nu\dt\bigl(g_{1-\alpha,n}*1\bigr)(t) + K_n(t).
	\end{aligned}
$$
Thus we obtain
\begin{equation}
\label{eq:4.10}
	\dt\bigl(g_{1-\alpha,n}*[\EE - (\EE_0+\nu)]\bigr)(t)
	+ \lambda t^{-\alpha\left(1-\frac{1}{m}\right)}\EE^m(t)
	\le L_n(t),
	\quad \mbox{a.e. $0<t<T$},
\end{equation}
where $L_n := I_n + J_n + K_n$ satisfying $L_n \to 0$ strongly in $L^1(0,T)$ as $n\to\infty$.
Let $Z \in C([0,\infty))$ be the solution to the Cauchy problem
$$
	\RL(Z-Z_0) + \lambda t^{-\alpha\left(1-\frac{1}{m}\right)}Z^m = 0, \quad t>0,
	\qquad Z(0) = Z_0 = \EE_0+\nu > 0.
$$
Convolving the equation for $Z$ with $h_{\alpha,n}$ defined in \eqref{eq:2.6} and taking the difference with \eqref{eq:4.10}, we obtain
\begin{equation}
\label{eq:4.11}
	\dt\bigl(g_{1-\alpha,n}*[\EE-Z]\bigr)(t)
	+ \lambda\left\{t^{-\alpha\left(1-\frac{1}{m}\right)}\EE^m(t) - \left(h_{\alpha,n}*\left[t^{-\alpha\left(1-\frac{1}{m}\right)}Z^m(\cdot)\right]\right)(t)\right\}
	\le L_n(t)
\end{equation}
for a.e. $0<t<T$.
To prove that $\EE(t) \le Z(t)$ for a.e. $0<t<T$, we assume that $\EE(t) > Z(t)$ for a.e. $0<t<T$.
We multiply \eqref{eq:4.11} by $[\EE-Z]_+ \in L^\infty(0,\infty)$
and apply Lemma~\ref{Lemma:2.1} for $H(y) = \frac{1}{2}[y]_+^2$ to obtain
$$
	\begin{aligned}
	&
	\frac{1}{2}\dt\bigl(g_{1-\alpha,n}*[\EE-Z]_+^2\bigr)(t) \\
	& \hspace{2cm}
	+ \lambda\left\{t^{-\alpha\left(1-\frac{1}{m}\right)}\EE^m(t) - \left(h_{\alpha,n}*\left[t^{-\alpha\left(1-\frac{1}{m}\right)}Z^m(\cdot)\right]\right)(t)\right\}[\EE(t)-Z(t)]_+ \\
	&
	\le L_n(t)[\EE(t)-Z(t)]_+, \quad \mbox{a.e. $0<t<T$}.
	\end{aligned}
$$
Since $\bigl(g_{1-\alpha,n}*[\EE-Z]_+^2\bigr)(0) = 0$,
by integrating both sides from $0$ to $T$, we have
\begin{equation}
\label{eq:4.12}
	\begin{aligned}
	&
	\frac{1}{2}\left(g_{1-\alpha,n}*[\EE-Z]_+^2\right)(T) \\
	& \hspace{1cm}
	+ \lambda\int_0^T\left\{t^{-\alpha\left(1-\frac{1}{m}\right)}\EE^m(t) - \left(h_{\alpha,n}*\left[t^{-\alpha\left(1-\frac{1}{m}\right)}Z^m(\cdot)\right]\right)(t)\right\}[\EE(t)-Z(t)]_+\,\dee t \\
	&
	\le \int_0^TL_n(t)[\EE(t)-Z(t)]_+\,\dee t.
	\end{aligned}
\end{equation}
Furthermore, by \eqref{eq:2.8} it holds that
$$
	\begin{aligned}
	&
	\left|\left(g_{1-\alpha,n}*[\EE-Z]_+^2\right)(T) - \left(g_{1-\alpha}*[\EE-Z]_+^2\right)(T)\right| \\
	&
	\le \|[\EE-Z]_+\|_{L^\infty(0,T)}^2\|g_{1-\alpha,n}-g_{1-\alpha}\|_{L^1(0,T)}
	\to 0 \quad \mbox{as $n\to\infty$}.
	\end{aligned}
$$
Since $Z \in L^\infty(0,\infty)$ and $t^{-\alpha\left(1-\frac{1}{m}\right)} \in L^1(0,T)$,
we see that $t^{-\alpha\left(1-\frac{1}{m}\right)}Z^m \in L^1(0,T)$, so it follows that
$$
	\begin{aligned}
	&
	\biggl|\int_0^T\left(h_{\alpha,n}*\left[t^{-\alpha\left(1-\frac{1}{m}\right)}Z^m(\cdot)\right]\right)(t)[\EE(t)-Z(t)]_+\,\dee t
	- \int_0^T t^{-\alpha\left(1-\frac{1}{m}\right)}Z^m(t)[\EE(t)-Z(t)]_+\,\dee t\biggr| \\
	&
	\le \|[\EE-Z]_+\|_{L^\infty(0,T)}\left\|\left(h_{\alpha,n}*\left[t^{-\alpha\left(1-\frac{1}{m}\right)}Z^m\right]\right) - t^{-\alpha\left(1-\frac{1}{m}\right)}Z^m\right\|_{L^1(0,T)} \\
	&
	\to 0 \quad \mbox{as $n\to\infty$}.
	\end{aligned}
$$
Moreover, since $L_n \to 0$ strongly in $L^1(0,T)$ as $n\to\infty$,
we also deduce that
$$
	\left|\int_0^TL_n(t)[\EE(t)-Z(t)]_+\,\dee t\right|
	\le \|[\EE-Z]_+\|_{L^\infty(0,T)}\|L_n\|_{L^1(0,T)}
	\to 0 \quad \mbox{as $n\to\infty$}.
$$
Thus, passing to the limit $n\to\infty$ in \eqref{eq:4.12},
we obtain
$$
	\frac{1}{2}\bigl(g_{1-\alpha}*[\EE-Z]_+^2\bigr)(T)
	+ \lambda\int_0^Tt^{-\alpha\left(1-\frac{1}{m}\right)}\bigl(\EE^m(t) - Z^m(t)\bigr)[\EE(t)-Z(t)]_+\,\dee t
	\le 0.
$$
Since $\bigl(\EE^m - Z^m\bigr)[\EE-Z]_+ \ge 0$ and $g_{1-\alpha}$ is nonnegative,
by arbitrariness of $T>0$, we see that $[\EE(t)-Z(t)]_+ = 0$ for a.e. $t>0$, that is,
$$
	\EE(t) \le Z(t), \quad \mbox{a.e. $t>0$}.
$$
Therefore, applying Proposition~\ref{Proposition:3.1} for $\beta = -\alpha\left(1-1/m\right) > -\alpha$ and $\gamma=m$,
we have
$$
	\EE(t) \le \frac{C}{1+t^{\frac{\alpha}{m^2}}}, \quad t>0.
$$
Furthermore, by assertion (i) of Lemma~\ref{Lemma:2.3}, we see that
$$
	\begin{aligned}
	\EE(t)
	&
	= \int_{\Omega}\left\{\frac{v^{m+1}(x,t)}{m+1} - \frac{S^{m+1}(x)}{m+1} - S^m(x)(v(x,t)-S(x))\right\}\,\dee x \\
	&
	\ge \frac{1}{m+1}\int_{\Omega}|v(x,t)-S(x)|^{m+1}\,\dee x \\
	&
	= \frac{1}{m+1}\|v(t)-S\|_{m+1}^{m+1},
	\quad \mbox{a.e. $t>0$}.
	\end{aligned}
$$
This implies that
$$
	\|v(t)-S\|_{m+1}^{m+1} \le \frac{C}{1+t^{\frac{\alpha}{m^2}}}, \quad t>0.
$$
For $1 \le p < m+1$, using H\"older's inequality and the fact that $\Omega$ is bounded, we deduce that
$$
	\|v(t)-S\|_p
	\le |\Omega|^{\frac{1}{p}-\frac{1}{m+1}}\|v(t)-S\|_{m+1}
	\to 0
	\quad \mbox{as $t\to\infty$}.
$$
Moreover, for $m+1 < p < \infty$, since $v$ and $S$ are bounded, we get
$$
	\|v(t)-S\|_p
	\le \|v(t)-S\|_\infty^{1-\frac{m+1}{p}}\|v(t)-S\|_{m+1}^{\frac{m+1}{p}}
	\le C\|v(t)-S\|_{m+1}^{\frac{m+1}{p}} \to 0
	\quad \mbox{as $t\to\infty$}.
$$
This implies that \eqref{eq:1.4} holds for $m\ge1$.
On the other hand,  by the fact that
$$
	|a^m-b^m| \le C_m(a^{m-1}+b^{m-1})|a-b|
$$
for $a,b \ge 0$ and $m\ge1$,
we also obtain
$$
	\|v^m(t)-S^m\|_p
	\le C\|v(t)-S\|_p
	\to 0 \quad \mbox{as $t\to\infty$}.
$$
Therefore, by assertion (ii) of Remark~\ref{Remark:1.1}, we obtain the convergence \eqref{eq:1.4} and \eqref{eq:1.5} for $m\ge1$.
$\Box$

\subsection{The fast diffusion case $0<m<1$}
\label{Subsection:4.3}
Next, we prove Theorem~\ref{Theorem:1.1} for the fast diffusion case $0<m<1$.
We first show that
\begin{equation}
\label{eq:4.13}
	\lambda\EE(t) \le \|\nabla(v^m(t)-S^m)\|_2^2, \quad \mbox{a.e. $0<t<T$},
\end{equation}
where $\lambda>0$ depends only on $\alpha$, $m$, $N$, $\Omega$, and $\|u_0\|_{L^{p_0}(\Omega)}$.
Indeed, by assertion (ii) of Lemma~\ref{Lemma:2.3}
and the Poincar\'{e} inequality, we deduce that
$$
	\begin{aligned}
	\EE(t)
	&
	= \int_{\Omega}\left\{\frac{v^{m+1}(x,t)}{m+1}-\frac{S^{m+1}(x)}{m+1} - S^m(x)(v(x,t)-S(x))\right\}\,\dee x \\
	&
	\le C\int_{\Omega}|v^m(x,t)-S^m(x)|^2\,\dee x \\
	&
	\le C\|\nabla(v^m(t)-S^m)\|_2^2, \quad \mbox{a.e. $0<t<T$}.
	\end{aligned}
$$

\subsubsection{The case $\alpha<m<1$}
Since $\alpha<m<1$, combining \eqref{eq:4.1}, the first estimate in \eqref{eq:4.2}, and \eqref{eq:4.13},
we obtain
$$
	\dt\bigl(g_{1-\alpha,n}*[\EE-(\EE_0+\nu)]\bigr)(t)
	+ \lambda t^{\alpha\left(\frac{1}{m}-1\right)}\EE(t)
	\le L_n(t),
	\quad \mbox{a.e. $0<t<T$}.
$$
Let $W \in C([0,\infty))$ be the solution to the Cauchy problem
$$
	\RL(W-W_0) + \lambda t^{\alpha\left(\frac{1}{m}-1\right)}W = 0, \quad t>0,
	\qquad W(0)=W_0=\EE_0+\nu>0.
$$
Repeating the same argument as in the case $m\ge1$, we have
$$
	\EE(t) \le W(t), \quad \mbox{a.e. $t>0$}.
$$
By Proposition~\ref{Proposition:3.1} for $\beta=\alpha\left(1/m-1\right)$ and $\gamma=1$,
we get
\begin{equation}
\label{eq:4.14}
	\EE(t) \le \frac{C}{1+t^{\frac{\alpha}{m}}}, \quad t>0.
\end{equation}
Furthermore, since $m<1$, it follows from Taylor's theorem that
\begin{equation}
\label{eq:4.15}
	\EE(t)
	= \frac{m}{2}\int_{\Omega}\eta^{m-1}(x,t)|v(x,t)-S(x)|^2\,\dee x
	\ge \frac{mM^{m-1}}{2}\|v(t)-S\|_2^2,
	\quad \mbox{a.e. $t>0$},
\end{equation}
where
$$
	M = \max\left\{\esssup_{t>0}\|v(t)\|_\infty, \|S\|_\infty\right\}.
$$
This together with \eqref{eq:4.14} implies that
$$
	\|v(t)-S\|_2^2 \le \frac{C}{1+t^{\frac{\alpha}{m}}}, \quad t>0.
$$
Therefore, using again boundedness of $v$ and $S$ and the fact that $\Omega$ is bounded,
we get \eqref{eq:1.4} for the case $\alpha<m<1$.
Moreover, by the elementary inequality
$$
	|a^m-b^m| \le C_m|a-b|^m, \quad a,b \ge 0,\ \ 0<m<1,
$$
we also obtain
$$
	\begin{aligned}
	\|v^m(t)-S^m\|_p^p
	&
	\le C_m\int_{\Omega}|v(x,t)-S(x)|^{pm}\,\dee x \\
	&
	\le C_m\left(\int_{\Omega}\,\dee x\right)^{1-m}\left(\int_{\Omega}|v(x,t)-S(x)|^p\,\dee x\right)^m \\
	&
	= C\|v(t)-S\|_p^{pm} \\
	&
	\to 0 \quad \mbox{as $t\to\infty$}.
	\end{aligned}
$$
This implies that \eqref{eq:1.5} also holds for $\alpha<m<1$.
$\Box$

\subsubsection{The case $0<m<\alpha$}

We next consider the case $0<m<\alpha$, namely, $\alpha/m>1$.
By combining \eqref{eq:4.1}, the last estimate in \eqref{eq:4.2} with $\beta=\alpha(1/m-1)>1-\alpha$, and \eqref{eq:4.13},
we also obtain
$$
	\dt\bigl(g_{1-\alpha,n}*[\EE-(\EE_0+\nu)]\bigr)(t) + \lambda t^{\beta}\EE(t)
	\le \mu t^{\beta-1} + L_n(t)
$$
for a.e. $t>0$.
Let $\tilde Z \in C([0,\infty))$ be the nonnegative solution to the Cauchy problem
$$
	\RL(\tilde Z - \tilde Z_0) + \lambda t^{\beta}\tilde Z = \mu t^{\beta-1},
	\quad t>0,
	\qquad \tilde Z(0) = \tilde Z_0 = \EE_0 + \nu.
$$
Then, using regularization arguments and Lemma~\ref{Lemma:2.1} again,
we have
$$
	\EE(t) \le \tilde Z(t), \quad t>0.
$$
Thus, applying Lemma~\ref{Lemma:3.2} for $\beta=\alpha\left(1/m-1\right)$,
we get
$$
	\EE(t) \le Ct^{-1}
$$
for sufficiently large $t>0$.
This together with \eqref{eq:4.15} implies that
$$
	\|v(t)-S\|_2^2 \le Ct^{-1}
$$
for sufficiently large $t>0$.
Thus the convergences \eqref{eq:1.4} and \eqref{eq:1.5} also hold for $0<m<\alpha$.
$\Box$

\subsubsection{The critical case $m=\alpha$}

Finally, we consider the critical case $m=\alpha$.
Combining \eqref{eq:4.1}, the second estimate in \eqref{eq:4.2}, and \eqref{eq:4.13}, we have
$$
	\dt\bigl(g_{1-\alpha,n}*[\EE-(\EE_0+\nu)]\bigr)(t)
	+ \lambda t^{1-\alpha}\EE(t)
	\le \mu t^{-\alpha}\log(1+t) + L_n(t)
$$
for a.e. $0<t<T$.
Let $\tilde{W} \in C([0,\infty))$ be the nonnegative solution to the Cauchy problem
$$
	\RL(\tilde{W}-\tilde{W}_0) + \lambda t^{1-\alpha}\tilde{W} = \mu t^{-\alpha}\log(1+t),
	\quad t>0,
	\qquad \tilde{W}(0)=\tilde{W}_0=\EE_0+\nu.
$$
Then we obtain
$$
	\EE(t) \le \tilde W(t), \quad t>0.
$$
By Lemma~\ref{Lemma:3.3} and \eqref{eq:4.15}, we have
$$
	\|v(t)-S\|_2^2 \le C\frac{\log(1+t)}{t}
$$
for sufficiently large $t>0$.
This implies that \eqref{eq:1.4} and \eqref{eq:1.5} also hold for the critical case $m=\alpha$.
Thus the proof of Theorem~\ref{Theorem:1.1} is complete.
$\Box$

\section{An alternative proof in the linear case $m=1$}
\label{Section:5}
In this section, we explain assertion (iii) of Remark~\ref{Remark:1.1}.
In the linear case $m=1$, the unique solution to problem~\eqref{P} admits the representation
$$
	u(x,t) = \sum_{k=1}^\infty E_{\alpha,1}(-\lambda_k t^{\alpha})\widehat u_{0,k}\phi_k(x),
$$
where $(\lambda_k,\phi_k)$ satisfy that
$$
	-\Delta\phi_k = \lambda_k\phi_k \quad \mbox{in $\Omega$},
	\qquad \phi_k = 0 \quad \mbox{on $\partial\Omega$},
$$
and $0 < \lambda_1 \le \lambda_2 \le \cdots$ (see, e.g., \cite{L}).
Moreover, the Fourier coefficients $\widehat{u}_{0,k}$ are defined by
$$
	\widehat{u}_{0,k} = \int_{\Omega}\phi_k(x)u_0(x)\,\dee x, \quad k=1,2,\cdots.
$$
On the other hand, $S$ also admits the representation
$$
	S(x) = \frac{1}{\Gamma(1-\alpha)}(-\Delta)^{-1}u_0(x)
	= \frac{1}{\Gamma(1-\alpha)}\sum_{k=1}^\infty\lambda_k^{-1}\widehat u_{0,k}\phi_k(x).
$$
Thus, since $\{\phi_k\}$ is an orthonormal basis of $L^2(\Omega)$, by Parseval's identity, we deduce that
$$
	\begin{aligned}
	\|t^{\alpha}u(t) - S\|_2^2
	&
	= \sum_{k=1}^\infty\left\|\widehat u_{0,k}\left(t^{\alpha}E_{\alpha,1}(-\lambda_k t^{\alpha}) - \frac{1}{\Gamma(1-\alpha)\lambda_k}\right)\phi_k\right\|_2^2 \\
	&
	= \sum_{k=1}^\infty|\widehat u_{0,k}|^2\left|t^{\alpha}E_{\alpha,1}(-\lambda_k t^{\alpha}) - \frac{1}{\Gamma(1-\alpha)\lambda_k}\right|^2.
	\end{aligned}
$$
Note that
$$
	E_{\alpha,1}(-\tau) \sim \frac{1}{\Gamma(1-\alpha)\tau} \quad \mbox{as $\tau\to\infty$}
$$
(see, e.g., \cite{BS}).
Hence, for each $k \in \N$,
$$
	\left|t^{\alpha}E_{\alpha,1}(-\lambda_k t^{\alpha}) - \frac{1}{\Gamma(1-\alpha)\lambda_k}\right|
	\to 0 \quad \mbox{as $t\to\infty$}.
$$
Furthermore, by \eqref{eq:2.1}, we obtain the following estimate uniformly in $t>0$
$$
	\begin{aligned}
	&
	\sum_{k=1}^\infty|\widehat u_{0,k}|^2\left|t^{\alpha}E_{\alpha,1}(-\lambda_k t^{\alpha}) - \frac{1}{\Gamma(1-\alpha)\lambda_k}\right|^2 \\
	&
	\le \sum_{k=1}^\infty|\widehat u_{0,k}|^2\left|\frac{C}{\lambda_k} - \frac{1}{\Gamma(1-\alpha)\lambda_k}\right|^2
	\le \sum_{k=1}^\infty|\widehat u_{0,k}|^2\frac{C}{\lambda_k^2} \\
	&
	\le \frac{C}{\lambda_1^2}\sum_{k=1}^\infty|\widehat u_{0,k}|^2
	\le \frac{C}{\lambda_1^2}\|u_0\|_2^2,
	\end{aligned}
$$
where we used Bessel's inequality in the last inequality.
Thus, by the dominated convergence theorem for series, we have
$$
	\|t^{\alpha}u(t)-S\|_2^2 \to 0
	\quad \mbox{as $t\to\infty$}.
$$
This completes the proof in the linear case.
$\Box$
\bigskip

\noindent
{\bf Acknowledgements.}
The author would like to express his sincere gratitude to Professor Tatsuki Kawakami for his helpful comments on this study.


\end{document}